\documentclass[10 pt, journal]{IEEEtran}  
\IEEEoverridecommandlockouts                             

\usepackage{xspace,amssymb,epsfig,subfigure,syntonly}
\usepackage{epsfig,amsmath,color}
\usepackage{xspace,syntonly}
\usepackage{url}
\usepackage{algorithm,algorithmic}
\usepackage{mathtools}
\usepackage{cite}
\usepackage{pstool}

\newcommand{\utwi}[1]{\mbox{\boldmath $#1$}}

\newcommand{\diag}{{\textrm{diag}}}

\newcommand{\cD}{{\cal D}}

\newcommand{\cL}{{\cal{L}}}
\newcommand{\cN}{{\cal N}}

\newcommand{\cG}{{\cal G}}
\newcommand{\cA}{{\cal A}}
\newcommand{\cS}{{\cal S}}

\newcommand{\cC}{{\cal C}}
\newcommand{\cE}{{\cal E}}

\newcommand{\cF}{{\cal F}}
\newcommand{\cB}{{\cal B}}

\newcommand{\cM}{{\cal M}}

\newcommand{\cY}{{\cal Y}}

\newcommand{\bc}{{\bf c}}
\newcommand{\ba}{{\bf a}}
\newcommand{\bb}{{\bf b}}

\newcommand{\be}{{\bf e}}
\newcommand{\bbf}{{\bf f}}
\newcommand{\bg}{{\bf g}}

\newcommand{\bp}{{\bf p}}
\newcommand{\bq}{{\bf q}}
\newcommand{\br}{{\bf r}}
\newcommand{\bs}{{\bf s}}
\newcommand{\bx}{{\bf x}}
\newcommand{\bu}{{\bf u}}

\newcommand{\bv}{{\bf v}}
\newcommand{\bi}{{\bf i}}
\newcommand{\bz}{{\bf z}}
\newcommand{\by}{{\bf y}}

\newcommand{\bB}{{\bf B}}

\newcommand{\bJ}{{\bf J}}
\newcommand{\bH}{{\bf H}}

\newcommand{\bR}{{\bf R}}

\newcommand{\bX}{{\bf X}}
\newcommand{\bZ}{{\bf Z}}

\newcommand{\bY}{{\bf Y}}

\newcommand{\bgamma}{{\utwi{\gamma}}}

\newcommand{\brho}{{\utwi{\rho}}}
\newcommand{\bupsilon}{{\utwi{\upsilon}}}

\newcommand{\bPhi}{{\utwi{\Phi}}}

\newcommand{\bxi}{{\utwi{\xi}}}

\newcommand{\bmu}{{\utwi{\mu}}}

\newcommand{\bvartheta}{{\utwi{\vartheta}}}
\newcommand{\bGamma}{{\utwi{\Gamma}}}

\newcommand{\sfH}{\textsf{H}}
\newcommand{\sfT}{\textsf{T}}

\usepackage{epstopdf}

\begin{document}

\newtheorem{definition}{Definition}
\newtheorem{remark}{Remark}
\newtheorem{proposition}{Proposition}
\newtheorem{lemma}{Lemma}
\newtheorem{theorem}{Theorem}
\def\HS{\hspace{\fontdimen2\font}}
\font\myfont=cmr12 at 16pt

\IEEEoverridecommandlockouts

\markboth{}%
{ \MakeLowercase{\textit{et al.}}: }


\title{Optimal Power Flow Pursuit}

\author{Emiliano Dall'Anese, \emph{Member, IEEE}, and Andrea Simonetto, \emph{Member, IEEE}
\thanks{IEEE Transactions on Smart Grid. Submitted: Jan. 25, 2016; revised: Apr. 8, 2016; accepted: May 15, 2016}
\thanks{The work of E. Dall'Anese and was supported in part by the Laboratory Directed Research and Development Program at the National Renewable Energy Laboratory.}
\thanks{E. Dall'Anese is  with the National Renewable Energy Laboratory, Golden, USA. A Simonetto is with the Universit\'e Catholique de Louvain, Louvain-la-Neuve, Belgium. E-mail: {emiliano.dallanese@nrel.gov}, {andrea.simonetto@uclouvain.be}.%
}
}
\maketitle

\maketitle

\begin{abstract}
This paper considers distribution networks featuring inverter-interfaced distributed energy resources, and develops distributed feedback controllers that continuously drive the inverter output powers to solutions of AC optimal power flow (OPF) problems. Particularly, the controllers update the power setpoints 
based on voltage measurements as well as given (time-varying) OPF targets, and entail elementary operations implementable onto low-cost microcontrollers that accompany power-electronics interfaces of gateways and inverters. The design of the control framework is based on suitable linear approximations of the AC power-flow equations as well as Lagrangian regularization methods. Convergence and OPF-target tracking capabilities of the  controllers are analytically established. Overall, the proposed method allows to bypass traditional hierarchical setups where feedback control and  optimization operate at distinct time scales, and to enable real-time optimization of distribution systems. 
\end{abstract}

\begin{IEEEkeywords} Distribution systems, optimal power flow, time-varying optimization, renewable integration, voltage regulation. 
\end{IEEEkeywords}

\section{Introduction}
\label{sec:Introduction}

The present paper seeks contributions in the domain of operation and control of distribution systems with high integration of  distributed energy resources. The objective is to develop distributed controllers that leverage the opportunities for fast feedback offered by power-electronics interfaced renewable energy sources (RESs), to continuously drive the system operation towards AC optimal power flow (OPF) targets. 
 
Prior works that focused on addressing power-quality and reliability concerns related to RES operating with business-as-usual practices~\cite{Liu08} have looked at the design of Volt/VAr, Volt/Watt, and droop-based control strategies to regulate output powers based on local measurements, so that terminal voltages are within acceptable levels (see, e.g.,~\cite{Aliprantis13, Tonkoski11,vonAppen14,Dorfler14}); these  strategies have the potential of controlling inverter outputs at a time scale that is consistent with the fast dynamics that govern the grid edge; however, they do not guarantee system-level optimality and stability claims are mainly based on empirical evidences.  On a different time scale, centralized and distributed OPF-type algorithms have been developed  for distribution systems to compute optimal steady-state inverter setpoints. Objectives of the OPF task at the distribution level include minimization of power losses as well as maximization of economic benefits to utility and end-users (e.g.,~\cite{Farivar12,Paudyal11,Robbins15,OID}); typical constraints in the OPF task ensure that voltage magnitudes and currents are within predetermined bounds, and RES setpoints are within given operational and hardware limits. It is well-known that the OPF problem is \emph{nonconvex} and NP-hard (see e.g.,~\cite{LavaeiLow}). Centralized approaches either utilize off-the-shelf solvers for nonlinear programs~\cite{Paudyal11,Khodr07}, or leverage convex relaxation and approximation techniques to obtain convex surrogates~\cite{Farivar12,OID,Robbins15,swaroop2015linear,LavaeiLow}. On the other hand, distributed solution approaches leverage the decomposability of the Lagrangian associated with convex reformulations/approximations of the OPF, and utilize iterative primal-dual-type methods to decompose the solution of the OPF task across devices~\cite{Robbins15,Erseghe14,Tse12}. 

OPF approaches have been successfully applied to optimize the operation of transmission systems. However, the  time required to collect all the problem inputs (e.g., loads across the network and available RES powers) and solve the OPF task may not be consistent with underlying distribution-systems dynamics. For example, Figure~\ref{Fig:F_load} provides a snapshot  of the loading of five secondary transformers located in a distribution feeder in Anatolia, CA~\cite{Bank13}; in this case, it is apparent that the inverter setpoints should be updated every second in order to cope with load variations and yet guarantee system-level optimality. However, existing distribution management systems (DMS) may not be able to solve the OPF task and dispatch setpoints in such a fast time scale. Distributed OPF approaches, where the power commands are updated at a slow time scale dictated by the convergence time of the distributed algorithm~\cite{Robbins15,Erseghe14,Tse12}, might systematically regulate the inverter power-outputs around outdated setpoints (possibly leading to violations of voltage and security limits).

In an effort to bypass traditional hierarchical setups where local feedback control and network optimization operate at distinct time scales~\cite{Dorfler14}, this paper develops a distributed control scheme that leverages the opportunities for fast feedback offered by power-electronics interfaced RESs, and continuously drives the inverter output powers towards OPF-based targets.  These targets capture well-defined performance objectives as well as voltage regulation constraints.  The design of the control framework is based on suitable linear approximations of the AC power-flow equations~\cite{swaroop2015linear} as well as the double-smoothing technique proposed in~\cite{Koshal11} for time-invariant optimization, and further extended to the time-varying setup in~\cite{SimonettoGlobalsip2014}. By virtue of this technical approach, the controllers entail elementary operations implementable into low-cost microcontrollers that accompany power-electronics interfaces of gateways and RESs. Further, while pursuing OPF solutions, the proposed controllers do not require knowledge of loads at all the feeder locations. Convergence and OPF-target tracking capabilities of the proposed controllers are analytically established. 

Prior efforts in this direction include e.g., the continuous-time feedback controllers that seek Karush-Kuhn-Tucker conditions for  economic dispatch optimality  for bulk power systems in~\cite{Jokic_JEPES}.  Recently, modified automatic generation and frequency control methods that incorporate optimization objectives corresponding to DC OPF problems have been proposed for lossless bulk power systems in e.g.,~\cite{NaLi_ACC14,Chen_CDC14}. A heuristic based on saddle-point-flow methods is utilized in~\cite{Elia-Allerton13} to synthesize controllers seeking AC OPF solutions. A droop-type control strategy for reactive power compensation in single-phase radial systems is proposed in~\cite{Zhang13} and convergence to a feasible power-flow solution is established; however, inverter capacity limits are not accounted for and loads are static. A local reactive power control strategy based on gradient-projection method is proposed in~\cite{Zhu15}, and convergence to the solution of a well-defined (static) optimization problem is studied. An online gradient algorithm for AC optimal power flow in single-phase radial networks is proposed in~\cite{LowOnlineOPF}; it is shown that the proposed algorithm converges to the set of local optima of a static AC OPF problem, and sufficient conditions under which the online OPF converges to a global optimum are provided.  A central controller for a number of resources in a feeder of microgrid is developed in~\cite{AndreayOnlineOpt}, based on continuous gradient steering algorithms; the framework accounts for errors in the implementable power setpoints, and convergence of the average setpoints to the minimum of the considered control objective is established. Finally, a reactive power control strategy is proposed in~\cite{Arnold15} for single-phase distribution systems with a tree topology based on an the so-called extremum-seeking control method. 

The proposed framework considerably broadens the approaches of~\cite{Jokic_JEPES,NaLi_ACC14, Elia-Allerton13,Zhang13,Zhu15,LowOnlineOPF} by focusing on AC OPF setups for distribution systems with arbitrary topologies and by establishing convergence and optimality in the case of time-varying loads and ambient conditions. The proposed approach offers significant contribution over~\cite{DhopleDKKT15} by establishing convergence results for the case of time-varying loads and ambient conditions and enabling low complexity implementations.     

The remainder of this paper is organized as follows. Section~\ref{sec:preliminariesandsystemmodel} outlines the system model and describes the target time-varying OPF problem. Section~\ref{sec:opfPursuit} addresses the synthesis of the proposed feedback controllers pursuing OPF solutions, and Section~\ref{sec:numericaltests} presents test cases. Finally, Section~\ref{sec:conclusions} concludes the paper. Relevant proofs are reported in the Appendix.

\section{Preliminaries and System model}
\label{sec:preliminariesandsystemmodel}

\subsection{System model}
\label{sec:systemmodel}

\begin{figure}[t] 
\hspace{-.4cm}\includegraphics[width=9.5cm]{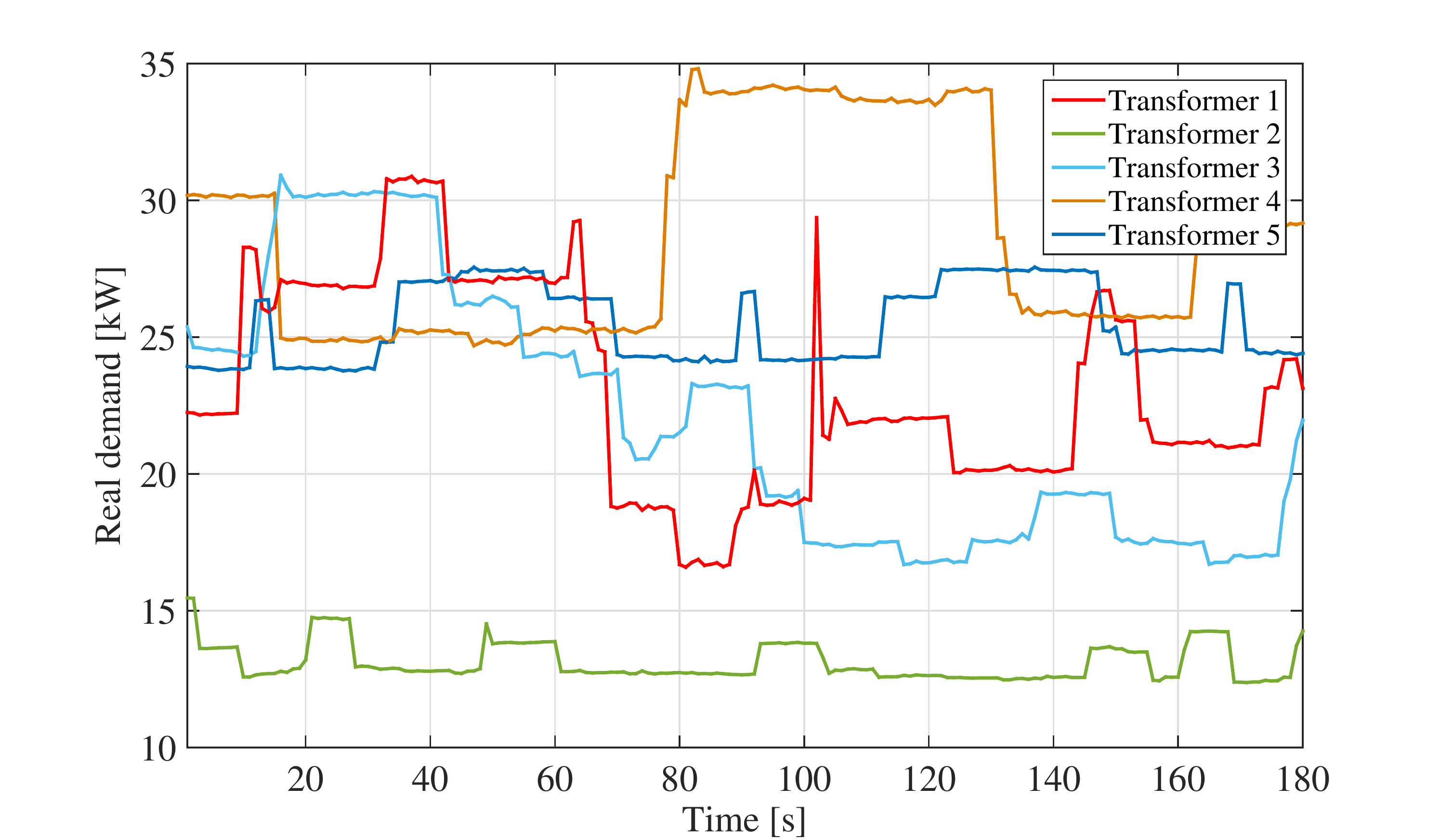}
\vspace{-.5cm}
\caption{Loading of five secondary transformers in a feeder located in Anatolia, CA~\cite{Bank13}. Trajectories correspond to transformer loadings from 5:00 PM to 5:03 PM during a day of August, 2012.}
\label{Fig:F_load}
\vspace{-.5cm}
\end{figure}

Consider a distribution feeder\footnote{Upper-case (lower-case) boldface letters will be used for matrices (column vectors); $(\cdot)^\sfT$ for transposition; $(\cdot)^*$
  complex-conjugate; and, $(\cdot)^\sfH$ complex-conjugate transposition; $\Re\{\cdot\}$ and $\Im\{\cdot\}$ denote the real and
  imaginary parts of a complex number, respectively; $\mathrm{j} := \sqrt{-1}$ the imaginary unit; and $|\cdot|$ denotes the absolute value of a number or the cardinality of a set. For $x \in \mathbb{R}$, function $[x]_+$ is defined as $[x]_+ := \max\{0,x\}$. For a given $N \times 1$ vector $\bx \in \mathbb{R}^N$, $\|\bx\|_2 := \sqrt{\bx^\sfH \bx}$; and, $\diag(\bx)$ returns a $N \times N$ matrix with the elements of $\bx$ in its diagonal. Further,  $\mathrm{proj}_{\cY}\{\bx\}$ denotes the projection of $\bx$ onto the convex set $\cY$. Given a given matrix $\bX \in \mathbb{R}^{N\times M}$, $x_{m,n}$  denotes its $(m,n)$-th entry. $\nabla_{\bx} f(\bx)$ returns the gradient vector of $f(\bx)$ with respect to $\bx \in \mathbb{R}^N$. Finally, $\mathbf{1}_N$ denotes the $N \times 1$ vector with all ones, and $\mathbf{0}_N$ denotes the $N \times 1$ vector with all zeros.} comprising $N+1$ nodes collected in the
set $\cN \cup \{0\}$, $\cN := \{1,\ldots,N\}$, and lines represented by the set of
edges $\cE := \{(m,n)\} \subset (\cN  \cup \{0\}) \times (\cN  \cup \{0\})$. Assume that  the temporal domain is discretized as $t = k \tau$, where $k \in \mathbb{N}$ and $\tau > 0$ is is a given interval, chosen to capture the variations on loads and ambient conditions [cf.~Figure~\ref{Fig:F_load}]. Let $V_n^k \in \mathbb{C}$ and $I_n^k \in \mathbb{C}$ denote the phasors for the line-to-ground voltage and the current injected at node $n$ over the $k$th slot, respectively, and define the $N$-dimensional complex vectors  $\bv^k := [V_1^k, \ldots, V_N^k]^\sfT \in \mathbb{C}^{N}$ and $\bi^k := [I_1^k, \ldots, I_N^k]^\sfT \in
\mathbb{C}^{N}$. Node $0$ denotes the secondary
of the distribution transformer, and it is taken to be the slack bus. Using Ohm's and Kirchhoff's circuit laws, the following linear relationship can be established:
\begin{align}
\left[
\begin{array}{c}
I_0^k \\
\bi^k
\end{array}
\right] = 
\underbrace{\left[
\begin{array}{cc}
 y_{00}^k  & (\bar{\by}^k)^\sfT  \\
 \bar{\by}^k  & \bY^k
\end{array}
\right]}_{:= \bY_{\mathrm{net}}^k}
\left[
\begin{array}{c}
V_0^k \\
\bv^k
\end{array}
\right]\ ,
\label{eq:iYv}
\end{align}
where the system admittance matrix $\bY_{\mathrm{net}}^k \in \mathbb{C}^{(N+1) \times (N+1)}$ is formed
based on the system topology and the $\pi$-equivalent circuit of the distribution lines (see e.g.,~\cite[Chapter 6]{kerstingbook} for additional details on distribution line modeling), and is partitioned in sub-matrices with  the following dimensions: $\bY^k \in \mathbb{C}^{N \times N}$, $\overline \by^k \in \mathbb{C}^{N \times 1}$, and $y_{00}^k \in \mathbb{C}$. Finally, $V_0^k = \rho_0 e^{\mathrm{j} \theta_0}$ is the slack-bus voltage with $\rho_0$ denoting the voltage magnitude at the secondary of the transformer.
A constant-power load model is utilized, and $P_{\ell,n}^k$ and $Q_{\ell,n}^k$ denote the real and reactive demands at node $n \in \cN$ at time $k$~\cite{kerstingbook}. 

Renewable energy sources (RESs) such as photovoltaic (PV) systems and small-scale wind turbines are assumed to be located at nodes $\cG \subseteq \cN$. For future developments, define $N_{\cG}:= |\cG|$. Given prevailing ambient conditions, let $P_{\textrm{av},n}^{k}$ denote the maximum renewable-based 
real power generation at node $n \in \cG$ at time $k$ -- hereafter referred to as the available real power. For example, for a PV system, the available power is a function of the incident irradiance, and corresponds to the maximum power point of the PV array. When RESs operate at unity power factor and inject in the network the whole available power,  a set of challenges related to power quality and reliability in distribution systems may emerge for sufficiently high levels of deployed RES capacity~\cite{Liu08}. For instance, overvoltages may be experienced during periods when RES generation exceeds the demand~\cite{Liu08}, while  fast-variations in the RES-output tend to propagate transients that lead to wear-out of legacy switchgear.  Efforts to ensure reliable operation of existing distribution systems with increased RES generation are focused on the possibility of inverters providing reactive power compensation and/or curtailing real power~\cite{Aliprantis13,Tonkoski11,Farivar12,OID}. Thus, letting $P_{n}^k$ and $Q_n^k$ denote the real and reactive powers at the AC side of inverter $n \in \cG$ at time $k$, the set of possible operating points can be specified as:
\begin{align}
 (P_n^k, Q_n^k) \in \cY_n^k & : =  \left\{({P}_{n}, {Q}_{n} ) \hspace{-.1cm} :  0 \leq {P}_{n}  \leq  P_{\textrm{av},n}^k, \right. \nonumber \\
&  \left.  \hspace{2cm} ({Q}_{n})^2  \leq  S_{n}^2 - ({P}_{n})^2 \right\} \label{mg-PV} 
 \end{align}
where $S_n$  is the rated apparent power. Lastly, the additional constraint $| {Q}_{n} | \leq (\tan \theta) {P}_{n}$ can be considered in the definition of $\cY_n^k$ to enforce a minimum power factor of $\cos \theta$; parameter $\theta$ can be conveniently tuned to account for a variety of control strategies, including reactive power compensation~\cite{Aliprantis13}, real power curtailment~\cite{Tonkoski11}, and joint real and reactive control~\cite{OID,DhopleDKKT15}. Other devices such as, for example, small-scale diesel generators, fuel cells, and variable speed drives can be  accommodated in the proposed framework by properly capturing their physical limits in the set $\cY_n^k$~\cite{DhopleNoFuel}. 

\subsection{Problem setup}
\label{sec:problem setup}

The objective is to develop distributed controllers that regulate the RES powers  $\{P_i^k, Q_i^k \}_{i \in \cG}$ at a time scale that is compatible with the distribution-systems dynamics, and operate in a closed-loop fashion as: 
\begin{subequations}
\label{eq:Objectivecontroller}
\begin{align}
[P_i^k, Q_i^k] & = \cC_i(P_i^{k-1}, Q_i^{k-1}, \by^k) , \hspace{.5cm} \forall i \in \cG  \label{eq:ObjectivecontrollerC} \\
\dot{\by}(t) & = \cF(\by, \{P_i^k, Q_i^k\}) \\
\by^k & = \cS(\by(t))\ ,
\end{align}
\end{subequations}
where $\cF(\cdot)$ models the physics of distribution systems (e.g., power flows) as well as the dynamics of primary-level inverter controllers~\cite{kerstingbook, Dorfler14,DhopleDKKT15}, $\by(t)$ represents pertinent electrical quantities (e.g., voltages and power flows), and $\by^k $ is a measurement of (some entries of) $\by(t)$ at time $k \tau$.  In the following, the control function $\cC_i(\cdot)$ will be designed in a way that the RES power outputs will continuously pursue solutions of an OPF problem. 

To this end, we begin with the formulation of a prototypical AC OPF problem, which is utilized to optimize the operation of the distribution feeder at time $k \tau$:
\begin{subequations} 
\label{Pmg}
\begin{align} 
 \mathrm{(OPF}^k \mathrm{)}  \hspace{1.6cm} & \hspace{-1.7cm} \min_{\bv, \bi,  \{P_i, Q_i \}_{i \in \cG} } \,\, h^k(\{V_i\}_{i \in \cN}) + \sum_{i \in \cG} f_i^k(P_i, Q_i)  \label{mg-cost} \\ 
& \hspace{-2.0cm} \mathrm{subject\,to}~\eqref{eq:iYv},  \mathrm{and}  \nonumber  \\ 
& \hspace{-2.8cm} V_i I_i^*  = P_i - P_{\ell,i}^k + \mathrm{j} (Q_i - Q_{\ell,i}^k), \hspace{.95cm}  \forall \, i \in \cG   \label{mg-balance-I} \\
& \hspace{-2.85cm} V_n I_n^* = - P_{\ell,n}^k - \mathrm{j} Q_{\ell,n}^k, \hspace{2.4cm}  \forall \, n \in \cN \backslash \cG \label{mg-balance-L} \\
& \hspace{-1.85cm} V^{\mathrm{min}} \leq |V_i| \leq V^{\mathrm{max}} ,  \hspace{1.7cm}  \forall \, i \in \cM  \hspace{-.2cm} \label{mg-Vlimits} \\
& \hspace{-1.4cm} (P_i, Q_i) \in \cY_i^k , \hspace{2.25cm} \forall \, i \in \cG \label{mg-PVp}\ ,
\end{align}
\end{subequations}
where $V^{\mathrm{min}}$ and $V^{\mathrm{max}}$ are minimum and maximum, respectively, voltage service limits (e.g., ANSI C.84.1 limits), $\cM \subseteq \cN$ is a set of nodes strategically selected to enforce voltage regulation throughout the feeder, $f_i^k(P_i, Q_i)$ is a time-varying function specifying performance objectives for the $i$th RES (e.g., cost of/reward for ancillary service provisioning~\cite{Farivar12,OID}, or feed-in tariffs~\cite{vonAppen14}), and $h^k(\{V_i\}_{i \in \cN})$ captures system-level objectives (e.g., power losses and/or deviations from the nominal voltage profile~\cite{OID}).  It is well-known that~\eqref{Pmg} is a \emph{nonconvex} (in fact, NP-hard) nonlinear program. Centralized~\cite{Farivar12,OID,Robbins15,swaroop2015linear,LavaeiLow} and \emph{distributed} solution approaches may not be able to solve $\mathrm{(OPF}^k \mathrm{)}$ and dispatch setpoints fast enough to cope with fast changes in the demand and ambient conditions at the grid edge (see e.g., Figure~\ref{Fig:F_load}), and might regulate the power-outputs $\{P_i, Q_i \}_{i \in \cG}$ around outdated setpoints (leading to suboptimal operation and potential violations of voltage and security limits). This is particularly relevant for distributed solution approaches, whereby the power commands are updated at a slow time scale, dictated by the convergence time of the distributed algorithm~\cite{Robbins15,Erseghe14,Tse12}. In contrast, the objective of~\eqref{eq:Objectivecontroller} is to update the power setpoints  at a fast time scale, and in a way that the inverter outputs are continuously regulated to a solution of $\mathrm{(OPF}^k \mathrm{)} $. How to design the control function~\eqref{eq:ObjectivecontrollerC} is the subject of the ensuing section.


\section{Design of feedback controllers}
\label{sec:opfPursuit}

\subsection{Leveraging approximate power-flow models}
\label{sec:approximate}

In this subsection, the linear approximation of the power flow equations proposed in~\cite{sairaj2015linear,swaroop2015linear} is briefly described; this approximation will be crucial to develop distributed feedback controllers that are low-complexity and fast acting. 

Let $\bs := [S_1, \ldots, S_N] \in \mathbb{C}^N$ collect the net power injected\footnote{For notational simplicity, in this subsection we drop the superscript $(\cdot)^k$ indexing the time instant $k \tau$ from all electrical and network quantities.} at nodes $\cN$, where $S_i = P_i - P_{\ell,i} + \mathrm{j} (Q_i - Q_{\ell,i})$ for $i \in \cG$, and $S_i = - P_{\ell,i} - \mathrm{j} Q_{\ell,i}$ for $i \in \cN \backslash  \cG$ [cf.~\eqref{mg-balance-I}--\eqref{mg-balance-L}]. Similarly, collect the voltage magnitudes $\{|V_i|\}_{i \in \cN}$ in $\brho :=  [|V_1|, \ldots, |V_N|]^\sfT \in \mathbb{R}^{N}$. The objective is to obtain  approximate power-flow relations whereby voltages  are \emph{linearly} related to injected powers $\bs$ as 
\begin{subequations} 
\label{eq:approximateVoltages}
\begin{align} 
\bv & \approx \bH \bp + \bJ \bq + \bb \label{eq:approximateV} \\
\brho & \approx \bR \bp + \bB \bq + \ba, \label{eq:approximate} 
\end{align}
\end{subequations}
where $\bp := \Re\{\bs\}$ and $\bq := \Im\{\bs\}$~\cite{sairaj2015linear,swaroop2015linear}. This way, the voltage constraints~\eqref{mg-Vlimits} can be approximated as $V^{\mathrm{min}} \mathbf{1}_N \leq \bR \bp + \bB \bq + \ba \leq V^{\mathrm{max}} \mathbf{1}_N$, while power-balance is intrinsically satisfied at all nodes; further, relevant electrical quantities of interests appearing in the function $h^k(\{V_i\}_{i \in \cN})$ in~\eqref{mg-cost}, e.g., power losses, can be expressed as linear functions of $\bp$ and $\bq$ (see e.g.,~\cite{swaroop2015linear}). What is more, by using~\eqref{eq:approximateV}--\eqref{eq:approximate}, function $h^k(\{V_i\})$ can be re-expressed as $\sum_{i\in\cG} h_i^k(P_i, Q_i)$. Following~\cite{sairaj2015linear,swaroop2015linear}, the matrices $\bR, \bB, \bH, \bJ$ and the vectors $\ba, \bb$ are obtained next.  

To this end, re-write~\eqref{mg-balance-I}--\eqref{mg-balance-L} in a compact form as 
\begin{align} 
\label{eq:s}
\hspace{-.1cm}
\bs =  \mathrm{diag}\left(\bv \right) \bi^* = \mathrm{diag}\left(\bv\right) (\bY^* \bv^* + \overline{\by}^* V_0^*) 
\end{align}
 and consider linearizing the AC power-flow equation around a given voltage profile $\bar{\bv} := [\bar{V}_1, \ldots, \bar{V}_N]^\sfT$~\cite{sairaj2015linear,swaroop2015linear}. In the following,  the voltages $\bv$ satisfying the nonlinear power-balance equations~\eqref{eq:s} are expressed as $\bv = \bar{\bv} + \be$, where the entries of $\be$ capture deviations around the linearization points $\bar{\bv}$. 
 For future developments, collect in the vector  $\bar{\brho} \in \mathbb{R}^N_+$ the magnitudes of voltages $\bar{\bv}$, and let $\bar{\bxi} \in \mathbb{R}^N$ and $\bar{\bvartheta} \in \mathbb{R}^N$ collect elements $\{\cos( \bar{\theta}_n)\}$ and $\{\sin( \bar{\theta}_n)\}$, respectively,    
 where $\bar{\theta}_i$ is the angle of the nominal voltage $\bar{V}_i$. 
 
By replacing $\bv$ with $\bar{\bv} + \be$ in~\eqref{eq:s} and discarding the second-order terms in $\be$ (e.g., discarding terms such as $\mathrm{diag}\left(\be\right) \bY^* \be^*$), equation~\eqref{eq:s} can be approximated as 
\begin{equation} 
\label{eq:Ve_approx}
\bGamma  \be + \bPhi \be^* = \bs + \bupsilon \, ,
\end{equation} 
where matrices $\bGamma$ and $\bPhi$ are given by $\bGamma := \mathrm{diag}\left(\bY^* \bar{\bv}^* + \overline{\by}^* V_0^* \right)$ and $\bPhi := \mathrm{diag}\left(\bar{\bv}\right) \bY^*$, respectively, and $\bupsilon := -\mathrm{diag}\left(\bar{\bv}\right) \left(\bY^* \bar{\bv}^*  + \overline{\by}^* V_0^*\right)$. Equation~\eqref{eq:Ve_approx} provides an approximate linear relationship between the injected complex powers and the voltage. In the following,~\eqref{eq:Ve_approx} will be further simplified by suitably selecting the nominal voltage profile $ \bar{\bv}$.  To this end, notice first that  matrix $\bY$ is invertible~\cite[Lemma~1]{sairaj2015linear}), and consider  the following choice of the nominal voltage $\bar{\bv}$:
\begin{equation} 
\label{eq:V_noload}
\bar{\bv} = - \bY^{-1} \overline{\by} V_0 \, .
\end{equation}
By using~\eqref{eq:V_noload}, one can see that $\bGamma = \mathbf{0}_{N \times N}$ and $\bupsilon = \mathbf{0}_N$, and therefore one obtains the linearized power-flow expression 
\begin{equation} \label{eq:linearized-power-flow}
\mathrm{diag}\left(\bar{\bv}^*\right) \bY \be = \bs^*. 
\end{equation}
A solution to~\eqref{eq:linearized-power-flow} can thus be expressed as $\be =  \bY^{-1}\mathrm{diag}^{-1}(\bar{\bv}^*) \bs^*$. Thus, expanding on this relation,  the approximate voltage-power relationship~\eqref{eq:approximateV} can be obtained by defining the matrices:   
\begin{subequations} 
\label{eq:Param_approx}
\begin{align}
& \hspace{-.2cm} \bar{\bR}  = \bZ_R \diag(\bar{\bxi}) (\diag(\bar{\brho}))^{-1} -\bZ_I \diag(\bar{\bvartheta}) (\diag(\bar{\brho}))^{-1} \\
& \hspace{-.2cm} \bar{\bB}  = \bZ_I \diag(\bar{\bxi}) (\diag(\bar{\brho}))^{-1} + \bZ_R \diag(\bar{\bvartheta}) (\diag(\bar{\brho}))^{-1} ,
\end{align}
\end{subequations}
where $\bZ_R := \Re\{\bY^{-1}\}$ and $\bZ_I := \Im\{\bY^{-1}\}$, and setting $\bH = \bar{\bR} + \mathrm{j} \bar{\bB}$, $\bJ = \bar{\bB} - \mathrm{j} \bar{\bR}$, and $\bb = \bar{\bv}$ . If the entries of $\bar{\bv}$ dominate those in $\be$, then $\bar{\brho} + \Re\{\be\}$ serves as a first-order approximation to the voltage magnitudes across the distribution network~\cite{sairaj2015linear}, and  relationship~\eqref{eq:approximateV} can be obtained by setting $\bR = \bar{\bR}$, $\bB = \bar{\bB}$, and $\ba = \bar{\brho}$. Analytical error bounds for~\eqref{eq:approximateV}--\eqref{eq:approximate} are provided in~\cite{sairaj2015linear};  the numerical experiments provided in~\cite{swaroop2015linear} demonstrate that~\eqref{eq:approximateV}--\eqref{eq:approximate} yield very accurate representations of the power flow equations. 



\subsection{Target time-varying optimization problem}
\label{sec:approximate}

To develop computationally affordable distributed controllers pursuing solutions to~\eqref{Pmg}, we begin with the derivation of a convex surrogate for the target OPF problem by leveraging~\eqref{eq:approximateVoltages} and~\eqref{eq:Param_approx}. Particularly, by using~\eqref{eq:approximate}, the voltage magnitude at node $n \in \cM$ and time $k$ can be approximated as $|V_n^k| \approx \sum_{i \in \cG} [r_{n,i}^k (P_i - P_{\ell,i}^k) + b_{n,i}^k (Q_i - Q_{\ell,i}^k)] + c_n^k$, with $c_n^k := \bar{\rho}_n^k - \sum_{i \in \cN \backslash \cG} (r_{n,i}^k P_{\ell,i}^k + b_{n,i}^k Q_{\ell,i}^k)$. It follows that problem~\eqref{Pmg} can be approximated as:
\begin{subequations} 
\label{Pmg2}
\begin{align} 
 \mathrm{(P1}^k \mathrm{)}  \hspace{1.8cm} & \hspace{-1.7cm} \min_{\{\bu_i \}_{i \in \cG} } \,\, \sum_{i \in \cG} \bar{f}_i^k(\bu_i) \label{mg-cost2} \\ 
& \hspace{-1.5cm} \mathrm{subject\,to}   \nonumber  \\ 
& \hspace{-0.8cm} g^k_n(\{\bu_i\}_{i \in \cG}) \leq 0 , \hspace{1.45cm} \forall n \in \cM \label{mg-volt1} \\
& \hspace{-0.8cm}  \bar{g}^k_n(\{\bu_i\}_{i \in \cG}) \leq 0 ,  \hspace{1.45cm} \forall n \in \cM \label{mg-volt2} \\
& \hspace{-0.8cm} \bu_i \in  \cY_i^k ,  \hspace{2.75cm} \forall \, i \in \cG\ , \label{mg-PVp2} 
\end{align}
\end{subequations}
where $\bu_i := [P_i, Q_i]^\sfT$, function $\bar{f}_i^k(\bu_i)$ is defined as $\bar{f}_i^k(\bu_i) := f_i^k(\bu_i) + h_i^k(\bu_i)$ for brevity, and
\begin{subequations} 
\label{eq:v} 
\begin{align} 
g^k_n(\{\bu_i\}_{i \in \cG}) & := V^{\mathrm{min}} - c_n^k \nonumber \\
&   - \sum_{i \in \cG} [r_{n,i}^k (P_i - P_{\ell,i}^k) + b_{n,i}^k (Q_i - Q_{\ell,i}^k)] \label{eq:vmin} \\
\bar{g}^k_n(\{\bu_i\}_{i \in \cG}) &:= \sum_{i \in \cG} [r_{n,i}^k (P_i - P_{\ell,i}^k) + b_{n,i}^k (Q_i - Q_{\ell,i}^k)] \nonumber \\
& + c_n^k  - V^{\mathrm{max}} \, . \label{eq:vmax} 
\end{align}
\end{subequations} 
Notice that the sets $\cY_i^k$, $i \in \cG$, are convex, closed, and bounded for all $k \geq 0$ [cf.~\eqref{mg-PV}]. For future developments, define the set $\cY^k := \cY_1^k \times \ldots \cY_{N_{\cG}}^k$. It is also worth reiterating that the $2 M$ constraints~\eqref{eq:v}, $M := |\cM|$, are utilized to enforce voltage regulation [cf.~\eqref{mg-Vlimits} and~\eqref{eq:approximate}].  Additional constraints can be considered in~$\mathrm{(OPF}^k \mathrm{)}$ and $\mathrm{(P1}^k \mathrm{)}$, but this would not affect the design of the feedback controllers. 

Regarding~\eqref{Pmg2}, the following  assumptions are made.   

\vspace{.1cm}

\noindent \emph{Assumption~1}.  Functions $f_i^k(\bu_i)$ and $h_i^k(\bu_i)$ are convex and continuously differentiable for each $i \in \cG$ and $k \geq 0$. Define further the gradient map:
\begin{align} 
\label{eq:gradientmap}
\bbf^k(\bu) := [\nabla_{\bu_1}^\sfT \bar{f}_1^k(\bu_1), \ldots, \nabla_{\bu_{N_{\cG}}}^\sfT \bar{f}_{N_{\cG}}^k(\bu_{N_{\cG}})]^\sfT \, .
\end{align}
Then, it is assumed that the gradient map $\bbf^k: \mathbb{R}^{2 N_{\cG}} \rightarrow \mathbb{R}^{2 N_{\cG}}$ is Lipschitz continuous
with constant $L$ over the compact set $\cY^k$ for all $k \geq 0$; that is, $\|\bbf^k(\bu) - \bbf^k(\bu^\prime) \|_2 \leq L \|\bu - \bu^\prime\|_2$, $\forall \,\, \bu, \bu^\prime \in \cY^k$. \hfill $\Box$ 

\vspace{.1cm}

\noindent \emph{Assumption~2 (Slater's condition)}. For all $k \geq 0$, there exist a set of feasible power injections $\{\hat{\bu}_i\}_{i \in \cG} \in \cY^{k}$ such that $g^k_n(\{\hat{\bu}_i\}_{i \in \cG}) \leq 0$ and $\bar{g}^k_n(\{\hat{\bu}_i\}_{i \in \cG}) \leq 0$, for all $n \in \cM$.  \hfill $\Box$ 

\vspace{.1cm}

Regarding \emph{Assumptions~2}, notice that functions $g^k_n(\{\hat{\bu}_i\}_{i \in \cG})$ and $\bar{g}^k_n(\{\hat{\bu}_i\}_{i \in \cG})$ are linear [cf.~\eqref{eq:v}]; hence, Slater's condition does not require strict inequalities~\cite{BoVa04}.  From the compactness of set $\cY^k$, and under \emph{Assumptions~1} and \emph{2}, problem~\eqref{Pmg2} is convex and strong duality holds \cite[Section~5.2.3]{BoVa04}. Further, there exists an optimizer at each time $k \geq 0$, which will be hereafter denoted as  $\{\bu_i^{\mathrm{opt},k}\}_{i \in \cG}$. For future developments, let $\bg^k(\bu) \in \mathbb{R}^M$ and $\bar{\bg}^k(\bu) \in \mathbb{R}^M$ be a vector stacking all the functions $g^k_n(\{\bu_i\}_{i \in \cG}), n \in \cM$, and $\bar{g}^k_n(\{\bu_i\}_{i \in \cG}), n \in \cM$, respectively; then, given that these functions are linear in $\bu$, it follows that there exists a constant $G$ such that $\|\nabla_\bu \bg^k(\bu)\|_2 \leq G$ and $\|\nabla_\bu \bar{\bg}^k(\bu)\|_2 \leq G$ for all $\bu \in \cY^k$ for all $k \geq 0$. 

It is worth pointing out that the cost functions $\{\bar{f}_i^k(\bu_i)\}_{i \in \cG}$ are not required to be strongly convex; whereas, the convergence properties of \emph{existing} distributed control schemes hinge on the strong convexity of the target cost functions (see e.g.,~\cite{Elia-Allerton13,NaLi_ACC14,Chen_CDC14}). 

Let $\cL^k(\bu, \bgamma, \bmu)$ denote the Lagrangian function associated with problem~\eqref{Pmg2}, where  $\bgamma := [\gamma_1, \ldots, \gamma_{M}]^\sfT$ and $\bmu := [\mu_1, \ldots, \mu_{M}]^\sfT$ collect the Lagrange multipliers associated with~\eqref{mg-volt1} and~\eqref{mg-volt2}, respectively. Further, let $\bu : = [(\bu_1)^\sfT, \ldots, (\bu_{N_{\cG}})^\sfT]^\sfT$ for brevity. Upon rearranging terms, the Lagrangian function can be expressed as 
\begin{align} 
& \cL^k(\bu, \bgamma, \bmu) := \sum_{i \in \cG} \bar{f}_i^k(P_i, Q_i) \nonumber \\
& + (P_i - P_{\ell,i}^k) (\check{\br}_i^k)^\sfT (\bmu - \bgamma) + (Q_i - Q_{\ell,i}^k) (\check{\bb}_i^k)^\sfT (\bmu - \bgamma) \nonumber\\ 
&  + (\bc^k)^\sfT(\bmu - \bgamma) + \bgamma^\sfT \mathbf{1}_m V^{\mathrm{min}}  - \bmu^\sfT \mathbf{1}_m V^{\mathrm{max}} 
\label{eq:lagrangianR]}
\end{align}
where $\check{\br}_i^k := [\{r_{j,i}^k\}_{j \in \cM}]^\sfT$ and $\check{\bb}_i^k := [\{b_{j,i}^k\}_{j \in \cM}]^\sfT$ are $M \times 1$ vectors collecting the entries of  $\bR^k$ and $\bB^k$ in the $i$th column and rows corresponding to nodes in $\cM$, and  $\bc^k := [\{c_{j}^k\}_{j \in \cM}]^\sfT$. Notice that, from the compactness of $\{\cY_i^k\}_{i \in \cG}$ and Slater's condition, it  follows that the optimal dual variables live in a compact set.

 In lieu of $\cL^k(\bu, \bgamma, \bmu)$, consider the following regularized Lagrangian function
\begin{align} 
& \cL_{\nu, \epsilon}^k(\bu, \bgamma, \bmu) := \cL^k(\bu, \bgamma, \bmu) \nonumber \\
& \hspace{2cm} + \frac{\nu}{2} \|\bu\|_2^2 - \frac{\epsilon}{2} (\|\bgamma\|_2^2 + \|\bmu\|_2^2) \label{eq:lagrangianR}
\end{align}
where the constant $\nu > 0$ and $\epsilon > 0$ appearing in the Tikhonov regularization terms are design parameters. Function~\eqref{eq:lagrangianR} is strictly convex  in the variables $\bu$ and strictly concave in the dual variables $\bgamma, \bmu$. The upshot of~\eqref{eq:lagrangianR} is that gradient-based approaches can be applied to~\eqref{eq:lagrangianR} to find an approximate solution to $\mathrm{(P1}^k \mathrm{)} $ with improved convergence properties~\cite{Koshal11,SimonettoGlobalsip2014}. Further, it allows one to drop the strong convexity assumption on $\{\bar{f}_i^k(\bu_i)\}_{i \in \cG}$ and to avoid averaging of primal and dual variables~\cite{OzdaglarSaddlePoint09}. Accordingly, consider the following saddle-point problem:
\begin{align} 
\label{eq:saddlepoint}
\max_{\bgamma \in \mathbb{R}^{M}_+, \bmu \in \mathbb{R}^{M}_+} \min_{\bu \in \cY^k}  \quad \cL_{\nu, \epsilon}^k(\bu, \bgamma, \bmu) 
\end{align}
and denote as $\{\bu_i^{*,k}\}_{i \in \cG}, \bgamma^{*,k}, \bmu^{*,k}$ the \emph{unique} primal-dual optimizer of~\eqref{eq:lagrangianR} at time $k$.  

In general, the solutions of~\eqref{Pmg2} and the regularized saddle-point problem~\eqref{eq:saddlepoint} are expected to be different; however, the discrepancy between $\bu_i^{\textrm{opt},k}$ and $\bu_i^{*,k}$ can be bounded as in~\cite[Lemma~3.2]{Koshal11}, whereas bounds of the constraint violation are substantiated in~\cite[Lemma~3.3]{Koshal11}. These bounds are proportional to $\sqrt{\epsilon}$; therefore, the smaller $\epsilon$, the smaller is the discrepancy between $\bu_i^{\textrm{opt},k}$ and $\bu_i^{*,k}$. 

Consider then the following primal-dual gradient method to solve the time-varying saddle-point problem~\eqref{eq:saddlepoint}:  
\begin{subequations} 
\label{eq:updateopt}
\begin{align} 
\tilde{\bu}_i^{k + 1} & =  \mathrm{proj}_{\cY_i^{k}}\left\{\tilde{\bu}_i^{k}  - \alpha \nabla_{\bu_i}  \cL_{\nu, \epsilon}^k(\bu, \bgamma, \bmu) |_{\tilde{\bu}_i^{k}, \tilde{\bgamma}^k, \tilde{\bmu}^k} \right\} , \nonumber \\
& \hspace{5cm} \forall \, i \in \cG \label{eq:updatev} \\
\tilde{\gamma}_n^{k + 1} & =   \mathrm{proj}_{\mathbb{R}_+} \left\{\tilde{\gamma}_n^{k} +  \alpha (g^k_n(\tilde{\bu}^k) - \epsilon \tilde{\gamma}_n^{k}) \right\} , \nonumber \\
& \hspace{5cm} \forall \, n \in \cM \label{eq:updategamma} \\
\tilde{\mu}_n^{k + 1} & =   \mathrm{proj}_{\mathbb{R}_+} \left\{\tilde{\mu}_n^{k} +  \alpha (\bar{g}^k_n(\tilde{\bu}^k) - \epsilon \tilde{\mu}_n^{k}) \right\}\ , \nonumber \\
& \hspace{5cm} \forall \, n \in \cM,  \label{eq:updatemu}
\end{align}
\end{subequations} 
where $\alpha > 0$ is the  stepsize and $\mathrm{proj}_{\cY}\{\bu\}$ denotes the projection of $\bu$ onto the convex set $\cY$; particularly, $\mathrm{proj}_{\mathbb{R}_+}\{u\} = \max\{0, u\}$, whereas~\eqref{eq:updatev} depends on the inverter operating region [cf.~\eqref{mg-PV}] and can be computed in closed-form (see e.g., Appendix~\ref{sec:setpointclosedform}). For the \emph{time-invariant} case (i.e., $\bar{f}_i^k(\bu_i) = \bar{f}_i(\bu_i)$, $g^k_n(\bu) = g_n(\bu)$, and $\bar{g}^k_n(\bu) = \bar{g}_n(\bu)$ for all $k >0$), convergence of~\eqref{eq:updateopt} is established in~\cite{Koshal11}. For the \emph{time-varying} case at hand, which captures the variability of underlying operating conditions at the grid edge [cf. Figure~\ref{Fig:F_load}], it is appropriate to  introduce additional assumptions to substantiate the discrepancy between the optimization problems that are associated with consecutive time instants~\cite{SimonettoGlobalsip2014}. 

\vspace{.1cm}

\noindent \emph{Assumption~3}. There exists a constant $\sigma_\bu \geq 0$ such that $\|\bu^{*,k+1} - \bu^{*,k}\| \leq \sigma_\bu$ for all $k \geq 0$. \hfill $\Box$ 

\vspace{.1cm}

\noindent \emph{Assumption~4}. There exist constants $\sigma_d \geq 0$ and $\sigma_{\bar{d}} \geq 0$ such that $|g^{k+1}_n(\bu^{*,k+1}) - g^{k}_n(\bu^{*,k}) | \leq \sigma_d$ and $|\bar{g}^{k+1}_n(\bu^{*,k+1}) - \bar{g}^{k}_n(\bu^{*,k}) | \leq \sigma_{\bar{d}}$, respectively, for all $n \in \cM$ and $k \geq 0$. \hfill $\Box$ 

\vspace{.1cm}

It can be shown that the conditions of \emph{Assumption~4} translate into bounds for the discrepancy between the optimal dual variables over two consecutive time instants; that is, $\|\bgamma^{*,k+1} - \bgamma^{*,k}\| \leq \sigma_{\bgamma}$ and $\|\bmu^{*,k+1} - \bmu^{*,k}\| \leq \sigma_{\bmu}$ with $\sigma_{\bgamma}$ and $\sigma_{\bmu}$ given by~\cite[Prop.~1]{SimonettoGlobalsip2014}. Upon defining $\bz^{*,k}:=[(\bu^{*,k})^\sfT, (\bgamma^{*,k})^\sfT, (\bmu^{*,k})^\sfT]^\sfT$  it also follows that $\|\bz^{*,k+1} - \bz^{*,k}\| \leq \sigma_{\bz}$ for a given $\sigma_{\bz} \geq 0$. Under \emph{Assumptions~1--4}, convergence of~\eqref{eq:updateopt} are investigated in~\cite[Theorem~1]{SimonettoGlobalsip2014}.

Similar to traditional distributed optimization schemes, updating the RES-power setpoints via~\eqref{eq:updateopt} leads to a setup where the optimization algorithm is \emph{decoupled} from the physical system~\cite{DhopleDKKT15}, and the RES-power setpoints are updated in an \emph{open-loop} fashion. In the next section, a feedback control architecture is proposed; actionable feedback from the distribution system will be incorporated in~\eqref{eq:updateopt} in order to enable adaptability to changing operating conditions.    

\begin{figure*}[t] 
\centering
\includegraphics[width=.75\textwidth, trim= 0cm .5cm 0cm .75cm, clip=on]{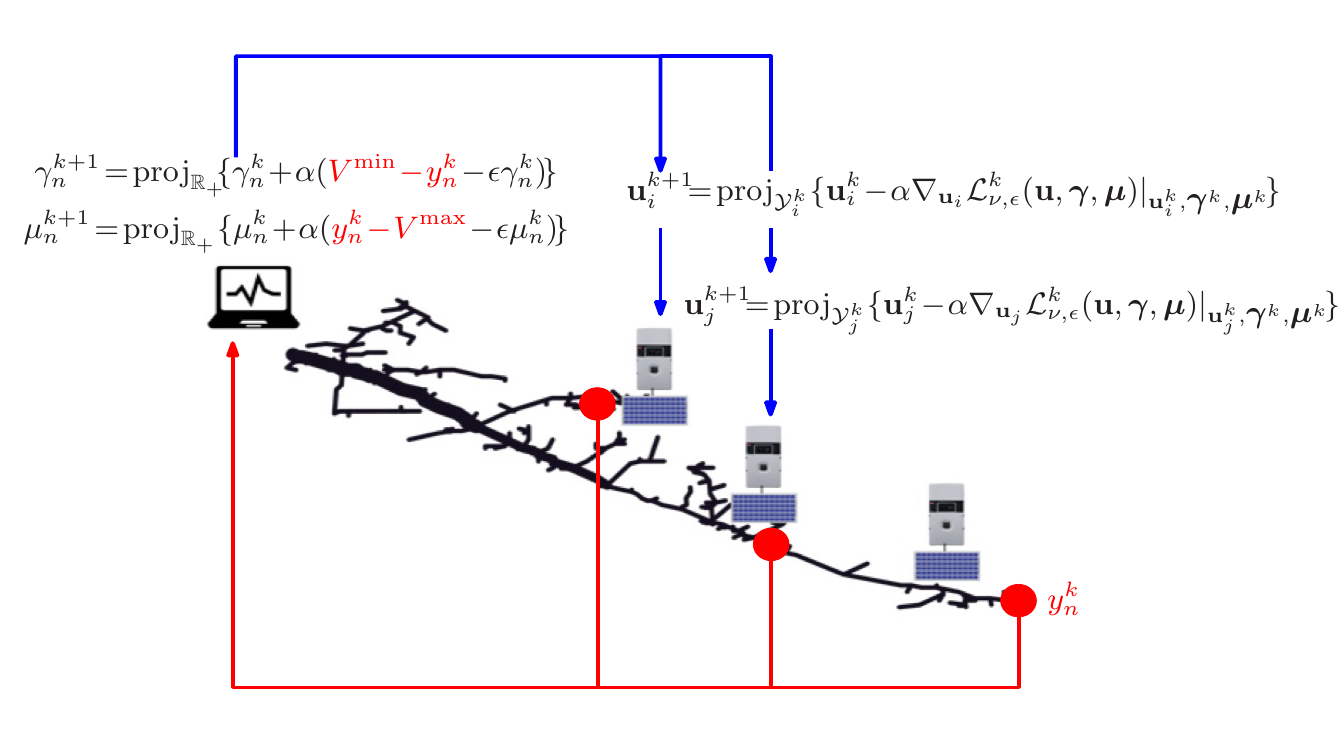}
\caption{Proposed control architecture. Upon collecting voltage measurements from selected feeder locations, RES-inverter setpoints are updated in a closed-loop fashion via~\eqref{eq:controllers}. }
\label{F_controller}
\end{figure*}

\subsection{Feedback controllers pursuing OPF solutions}
\label{sec:Pursuing}

Let $y_n^k$ denote a measurement of $|V_n^k|$ acquired at time $k$ from node $n \in \cM$ of the feeder. Then, we propose the following strategy to update the RES-inverter setpoints at each time $k$:

\vspace{.1cm}

\begin{subequations} 
\label{eq:controllers} 

\noindent \textbf{[S1]} Collect voltage measurements $\{y_n^k\}_{n\in\cM}$.

\noindent \textbf{[S2]} For all $n \in \cM$, update dual variables as follows:
\begin{align} 
\gamma_n^{k + 1} & =   \mathrm{proj}_{\mathbb{R}_+} \left\{\gamma_n^{k} +  \alpha (V^{\mathrm{min}} -  y_n^k - \epsilon \gamma_n^{k}) \right\}  \label{eq:updategammaC} \\
\mu_n^{k + 1} & =   \mathrm{proj}_{\mathbb{R}_+} \left\{\mu_n^{k} +  \alpha (y_n^k - V^{\mathrm{max}}  - \epsilon \mu_n^{k}) \right\}\ .  \label{eq:updatemuC}
\end{align}

\noindent \textbf{[S3]} Update power setpoints at each RES $i \in \cG$ as:   
\begin{align} 
\bu_i^{k + 1} & =  \mathrm{proj}_{\cY_i^{k}} \left\{\bu_i^{k}  - \alpha \nabla_{\bu_i}  \cL_{\nu, \epsilon}^k(\bu, \bgamma, \bmu) |_{\bu_i^{k}, \bgamma^k, \bmu^k} \right\} \label{eq:updatevC} 
\end{align}
\end{subequations}
and go to \textbf{[S1]}. 

\vspace{.1cm}

The proposed feedback control strategy is illustrated in Figure~\ref{F_controller}. It can be seen that the update~\eqref{eq:updatevC} is performed locally at each RES inverter and affords a closed-form solution for a variety of sets $\cY_i^k$ (see e.g., Appendix~\ref{sec:setpointclosedform}); updates~\eqref{eq:updategammaC}--\eqref{eq:updatemuC} can be computed either at each inverter (if the voltage measurements are broadcasted to the RES units) or at the utility/aggregator. The controller~\eqref{eq:updatevC} produces a  (continuous-time)  reference signal $\bu_i(t)$ for RES $i$ that is has step changes at instants $\{\tau k\}_{k \geq 0}$, is a left-continuous function, and takes the constant value $\bu_i^{k + 1}$ over the time interval $(\tau k, \tau (k+1)]$. Notice that differently from traditional distributed optimization schemes,~\eqref{eq:controllers} \emph{does not} require knowledge of the loads at locations $\cN \backslash \cG$. The only information required by the controllers pertains to the line and feeder models, which are utilized to build the network-related matrices in~\eqref{eq:approximateVoltages}. In the following, the convergence properties of the proposed scheme are analyzed. 

Key to this end is to notice that steps~\eqref{eq:updategammaC}--\eqref{eq:updatemuC} are in fact $\varepsilon$-gradients of the regularized Lagrangian function~\cite{Bertsekas1999}; that is, $V^{\mathrm{min}} -  y_n^k - \epsilon \gamma_n^{k} \neq \nabla_{\gamma_n}  \cL_{\nu, \epsilon}^k |_{\bu^{k}, \bgamma^k, \bmu^k}$ and $y_n^k - V^{\mathrm{max}}  - \epsilon \mu_n^{k} \neq \nabla_{\mu_n}  \cL_{\nu, \epsilon}^k|_{\bu^{k}, \bgamma^k, \bmu^k}$. This is primarily due to i) voltage measurements errors, ii) approximation errors introduced by~\eqref{eq:approximate}, and iii) setpoints possibly updated at a faster rate that the power-output settling time for off-the-shelf inverters~\cite{DhopleDKKT15}. The latter point is particularly important because updates~\eqref{eq:controllers} can be conceivably performed at a fast time scale (e.g., $\tau$ can be on the order of the subsecond); in fact, iterates $\bu_i^{k + 1}$, $\gamma_n^{k + 1}$, and $\mu_n^{k + 1}$ are updated via basic mathematical operations, and low latencies can be achieved with existing communications technologies. 

 Let $\be_\gamma^k \in \mathbb{R}^M$ and $\be_\mu^k \in \mathbb{R}^M$  collect the dual gradient \emph{errors} $V^{\mathrm{min}} -  y_n^k - \epsilon \gamma_n^{k} - \nabla_{\gamma_n}  \cL_{\nu, \epsilon}^k $ and $y_n^k - V^{\mathrm{max}}  - \epsilon \mu_n^{k} - \nabla_{\mu_n}  \cL_{\nu, \epsilon}^k$, respectively. Then, the following practical assumption is made.

\vspace{.1cm}

\noindent \emph{Assumption~5}. There exist a constant $e \geq 0$ such that $\max\{\|\be_\gamma^k\|_2, \|\be_\mu^k\|_2\}  \leq e$ for all $k \geq 0$. \hfill $\Box$ 

\vspace{.1cm}

Before stating the main convergence result for the network feedback controllers illustrated in Figure~\ref{F_controller}, it is convenient to introduce relevant definitions as well as a supporting lemma. Recall that $\tilde{\bz}^{k}:=[(\tilde{\bu}^{k})^\sfT, (\tilde{\bgamma}^{k})^\sfT, (\tilde{\bmu}^{k})^\sfT]^\sfT$, and define the time-varying mapping $\bPhi^k$ as
\begin{equation*}
\label{phi_mapping}
 \bPhi^k: \{\tilde{\bu}^{k}, \tilde{\bgamma}^{k}, \tilde{\bmu}^{k}\} \mapsto 
 \left[\begin{array}{c}
 \nabla_{\bu_1}  \cL_{\nu, \epsilon}^k(\bu, \bgamma, \bmu) |_{\tilde{\bu}_1^{k}, \tilde{\bgamma}^k, \tilde{\bmu}^k} \\
 \vdots \\
 \nabla_{\bu_{N_{\cG}}}  \cL_{\nu, \epsilon}^k(\bu, \bgamma, \bmu) |_{\tilde{\bu}_{N_{\cG}}^{k}, \tilde{\bgamma}^k, \tilde{\bmu}^k}\\
- (g^k_1(\tilde{\bu}^k) - \epsilon \tilde{\gamma}_1^{k})  \\
\vdots \\
- (g^k_M(\tilde{\bu}^k) - \epsilon \tilde{\gamma}_M^{k})  \\
 - (\bar{g}^k_1(\tilde{\bu}^k) - \epsilon \tilde{\mu}_1^{k}) \\
 \vdots\\
  - (\bar{g}^k_M(\tilde{\bu}^k) - \epsilon \tilde{\mu}_M^{k})  \end{array}\right] \hspace{-.1cm}
\end{equation*}
which is utilized to compute the gradients in the error-free iterates~\eqref{eq:updateopt} as
\begin{equation}\label{eq:updateopt_compact}
\tilde{\bz}^{k+1} = \mathrm{proj}_{\cY^{k} \times \mathbb{R}_{+}^M \times \mathbb{R}_{+}^M}\left\{\tilde{\bz}^k - \alpha \bPhi^k(\tilde{\bz}^k)\right\}.
\end{equation}
Then, the following holds.

\vspace{.1cm}

\begin{lemma}
\label{lemma-Phi}
The map $\bPhi^k$ is strongly monotone with constant $\eta = \min\{\nu,\epsilon\}$, and Lipschitz over $\cY^k \times \mathbb{R}^M_+ \times \mathbb{R}^M_+$ with constant $L_{\nu, \epsilon} = \sqrt{(L+\nu+2G)^2 + 2(G+\epsilon)^2}$. \hfill $\Box$
\end{lemma}

\vspace{.1cm}

The result above is a relaxed version of~\cite[Lemma~3.4]{Koshal11}, since it does not require the Lipschitz continuity of the gradient of~\eqref{mg-volt1}--\eqref{mg-volt2}. Convergence and tracking properties of the feedback controllers~\eqref{eq:controllers} are established next.

\vspace{.1cm}

\begin{theorem}
\label{theorem.inexact}
Consider the sequence $\{\bz^k\}: = \{\bu^k, \bgamma^k, \bmu^k\}$ generated by~\eqref{eq:controllers}. Let \emph{Assumptions~1--5} hold. For fixed  positive scalars $\epsilon, \nu >0$, if the stepsize $\alpha>0$ is chosen such that
\begin{equation}\label{eq.alpha}
\rho(\alpha):= \sqrt{1 - 2 \eta  \alpha + \alpha^2 L_{\nu, \epsilon}^2} < 1,
\end{equation}
that is $0 <\alpha < 2 \eta/ L_{\nu, \epsilon}^2$, then the sequence $\{\bz^k\}$ converges Q-linearly to $\bz^{*,k} := \{\bu^{*,k}, \bgamma^{*,k}, \bmu^{*,k}\}$ up to the asymptotic error bound given by: 
\begin{align}
\limsup_{k\to \infty} \|\bz^k - \bz^{*,k}\|_2 = \frac{1}{1 - \rho(\alpha)} \Big[\sqrt{2}\alpha e + \sigma_{\bz}\Big]. 
\label{eq.asympt_error}
\end{align}

\emph{Proof}. See the Appendix. \hfill $\Box$
\end{theorem}

\vspace{.1cm}

Equation~\eqref{eq.asympt_error} quantifies the maximum discrepancy between the iterates $ \{\bu^k, \bgamma^k, \bmu^k\}$ generated by the proposed controllers and the (time-varying) optimizer of problem~\eqref{eq:saddlepoint}. From~\cite[Lemma~3.2]{Koshal11} and by using the triangle inequality, a bound for the difference between $\bu^k$ and the time-varying solution of~\eqref{Pmg2} can be obtained.   The condition~\eqref{eq.alpha} imposes the requirements on the stepsize $\alpha$, such that $\rho(\alpha)$ is strictly less than $1$ and thereby enforcing Q-linear convergence. The optimal stepsize selection for convergence is $\alpha = \eta/L_{\nu, \epsilon}^2$.

The error~\eqref{eq.asympt_error} provides trade-offs between smaller $\alpha$'s (leading to a smaller term multiplying the gradient error $e$, and yet yielding poorer convergence properties, i.e., $\rho(\alpha)$ close to $1$) and bigger $\alpha$'s (leading to the opposite).

\vspace{.1cm}

\noindent \emph{Remark~1}. For notational and exposition simplicity, the paper  considered a balanced distribution network. However, the proposed control framework is applicable to multi-phase unbalanced systems with any topology. In fact, the linearized model in Section~\ref{sec:approximate} can be readily extended to the multi-phase unbalanced setup, and the controllers~\eqref{eq:controllers} can be embedded into inverters located at any phase and node. \hfill $\Box$

\vspace{.1cm}

\noindent \emph{Remark~2}. \emph{Assumption~2} requires the objective function~\eqref{mg-cost2} to be continuously differentiable.  Notice however that non-differentiable functions such as $|x|$ and $[x]_+ := \max\{0,x\}$ (with the latter playing an important role when feed-in tariffs are considered~\cite{vonAppen14}) can be readily handled upon introducing auxiliary optimization variables along with appropriate inequality constraints. For example, the problem $\min_x [x]_+$ s.t. $g(x) \leq 0$ can be reformulated in the following equivalent way:  $\min_{x, z} z$ s.t. $g(x) \leq 0, x \leq z$, and $z \geq 0$. \hfill $\Box$

\vspace{.1cm}

\noindent \emph{Remark~3}. Traditional OPF approaches include voltage regulation constraints at all nodes\cite{Farivar12,OID,Robbins15,swaroop2015linear,LavaeiLow}. In the present setup, the set $\cM$ corresponds to $M$ nodes where voltage measurements can be collected and utilized as actionable feedback in~\eqref{eq:controllers}. Accordingly, the set $\cM$ may include: i)  nodes $\mathcal{G}$ where RESs are located (existing inverters that accompany RESs are equipped with modules that measure the voltage at the point of connection); and, ii) additional nodes of a distribution feeder where distribution system operators deploy communications-enabled meters for voltage monitoring. 
\hfill $\Box$

\vspace{.1cm}

\noindent \emph{Remark~4}. The scalars $\sigma_\bu$, $\sigma_d$ and $\sigma_{\bar{d}}$ (and, thus, $\sigma_\bz$) in \emph{Assumption~3} and \emph{Assumption~4}  quantify the variability of the ambient and network conditions over the time interval $[\tau k, \tau (k+1)]$ as well as the (maximum) discrepancy between OPF solutions corresponding to two consecutive time instants $\tau k$ and $\tau (k+1)$ [cf. Figure~\ref{F_controller}]. On the other hand, parameter $e$ implicitly bounds the error between the setpoint $\bu_i^k$ commanded to the inverter and the actual inverter output, and it is related to the inverter's actuation time. It is worth pointing out that the results of Theorem~\ref{theorem.inexact} hold for \emph{any} value of $\tau$ (and, hence, for any values of the scalars $\sigma_\bu$, $\sigma_d$, $\sigma_{\bar{d}}$, and $\sigma_\bz$) and for any value of $e$. For given dynamics of ambient conditions, network, and problem parameters, $\sigma_\bz$ and $e$ are utilized to characterize the performance of the proposed controllers. For example, it is clear that the value of $\sigma_\bu$ decreases with $\tau$; and, as a consequence,  the distance between the controller output and the OPF solutions decreases with $\tau$ too [cf.~\eqref{eq.asympt_error}].  \hfill $\Box$

\section{Example of application}
\label{sec:numericaltests}

As an  application, a distribution network with high-penetration of photovoltaic (PV) systems is considered; particularly, it is  demonstrated how the proposed controllers can reliably prevent overvoltages that are likely to be experienced during periods when PV generation exceeds the demand~\cite{Liu08}.  

\begin{figure}[t] 
\centering
%
\vspace{.25cm}
\includegraphics[width=.45\textwidth]{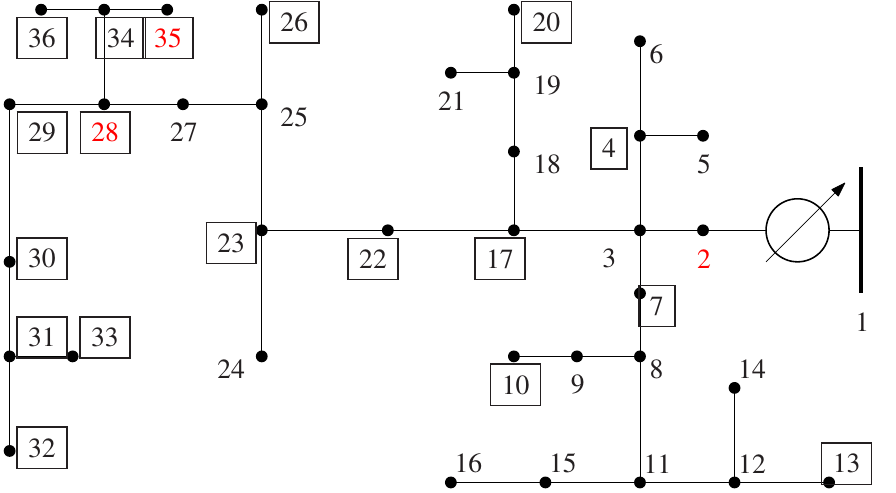}
\caption{IEEE 37-node feeder. The boxed nodes represent the location of PV systems. The red nodes are the ones analyzed in the numerical example.}
\label{F_feeder}
\end{figure}

\begin{figure}[t] 
\subfigure[]{\hspace{-.2cm}\includegraphics[width=9.0cm]{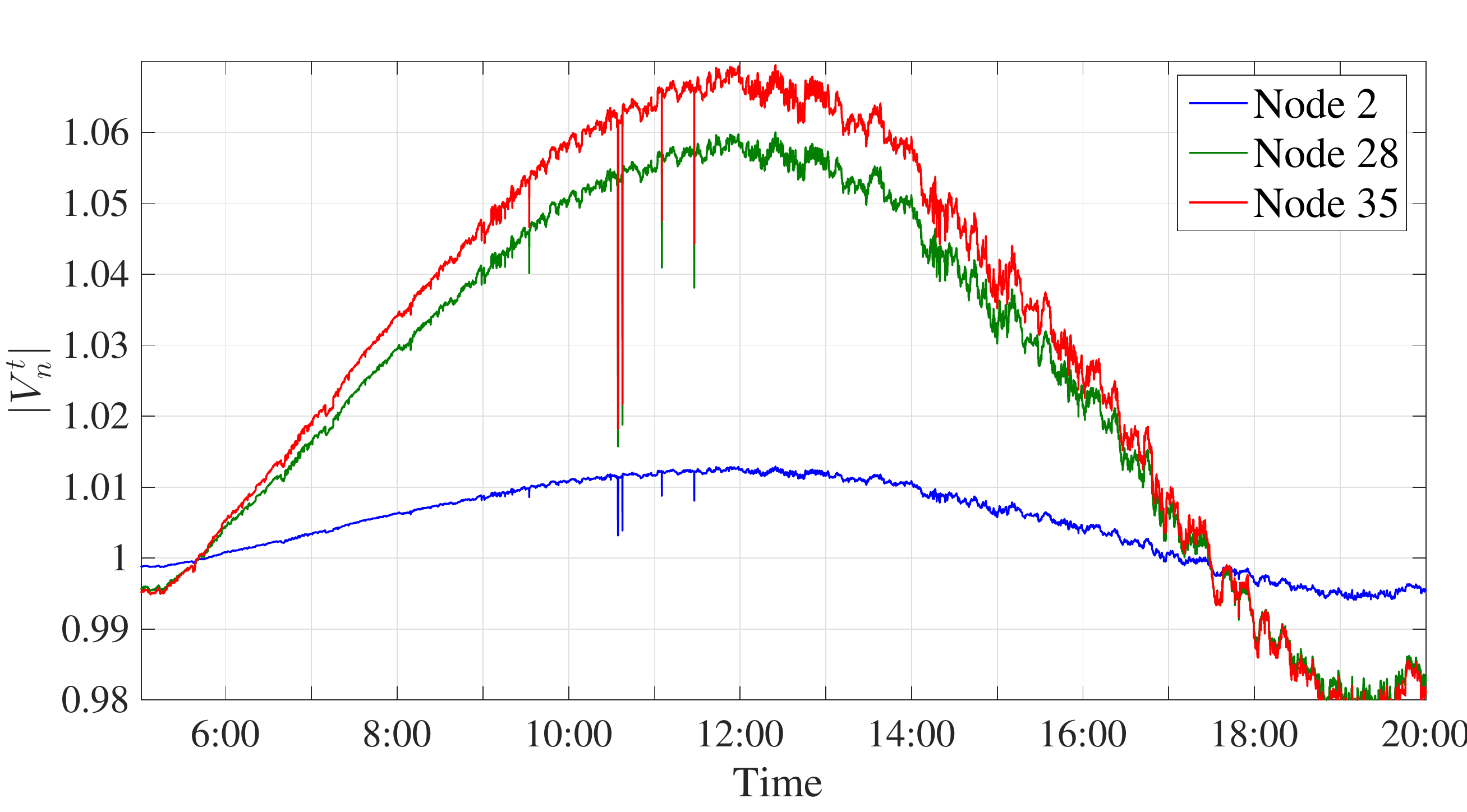} }
\subfigure[]{\hspace{-.2cm}\includegraphics[width=9.0cm]{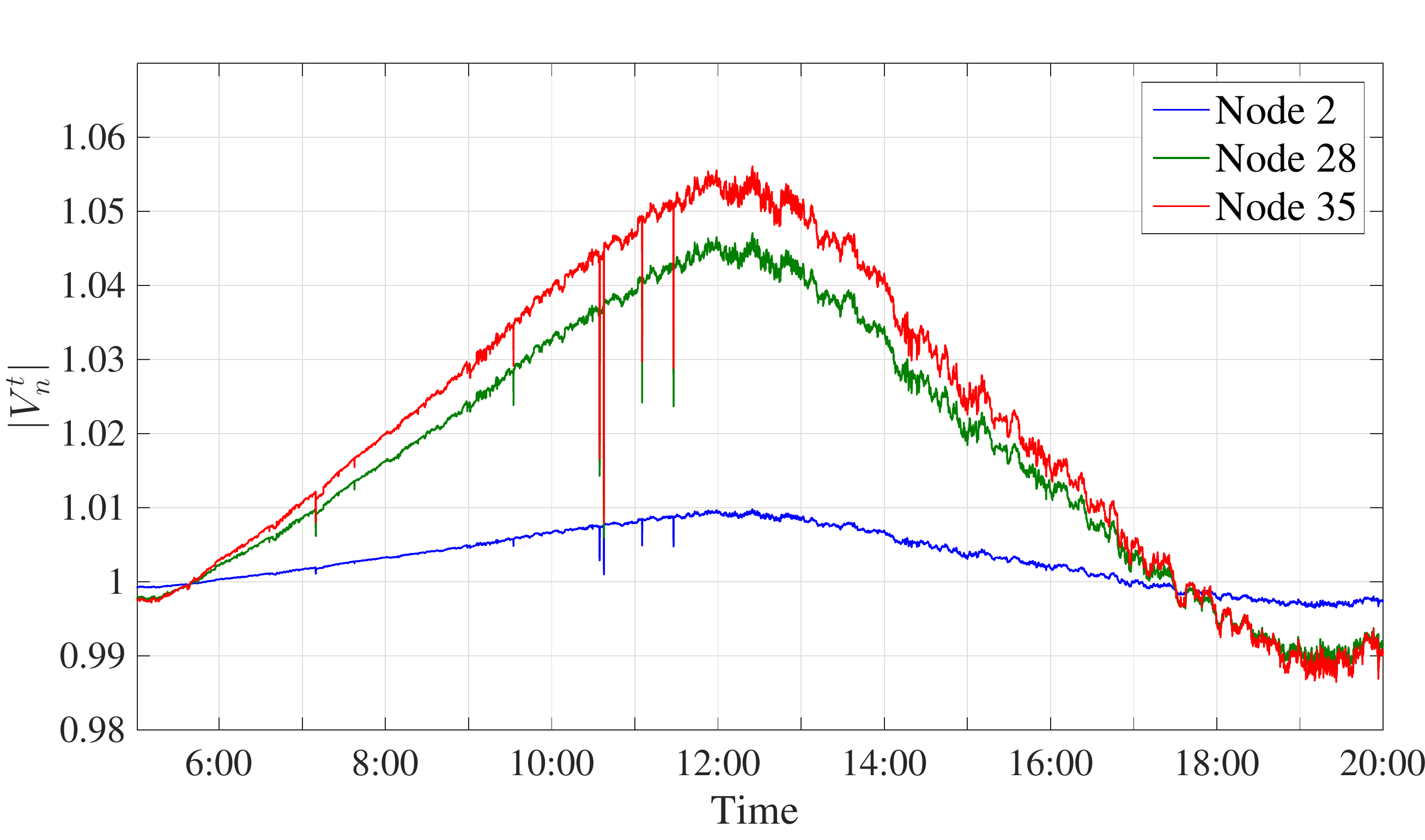} }
\subfigure[]{\hspace{-.2cm}\includegraphics[width=9.0cm]{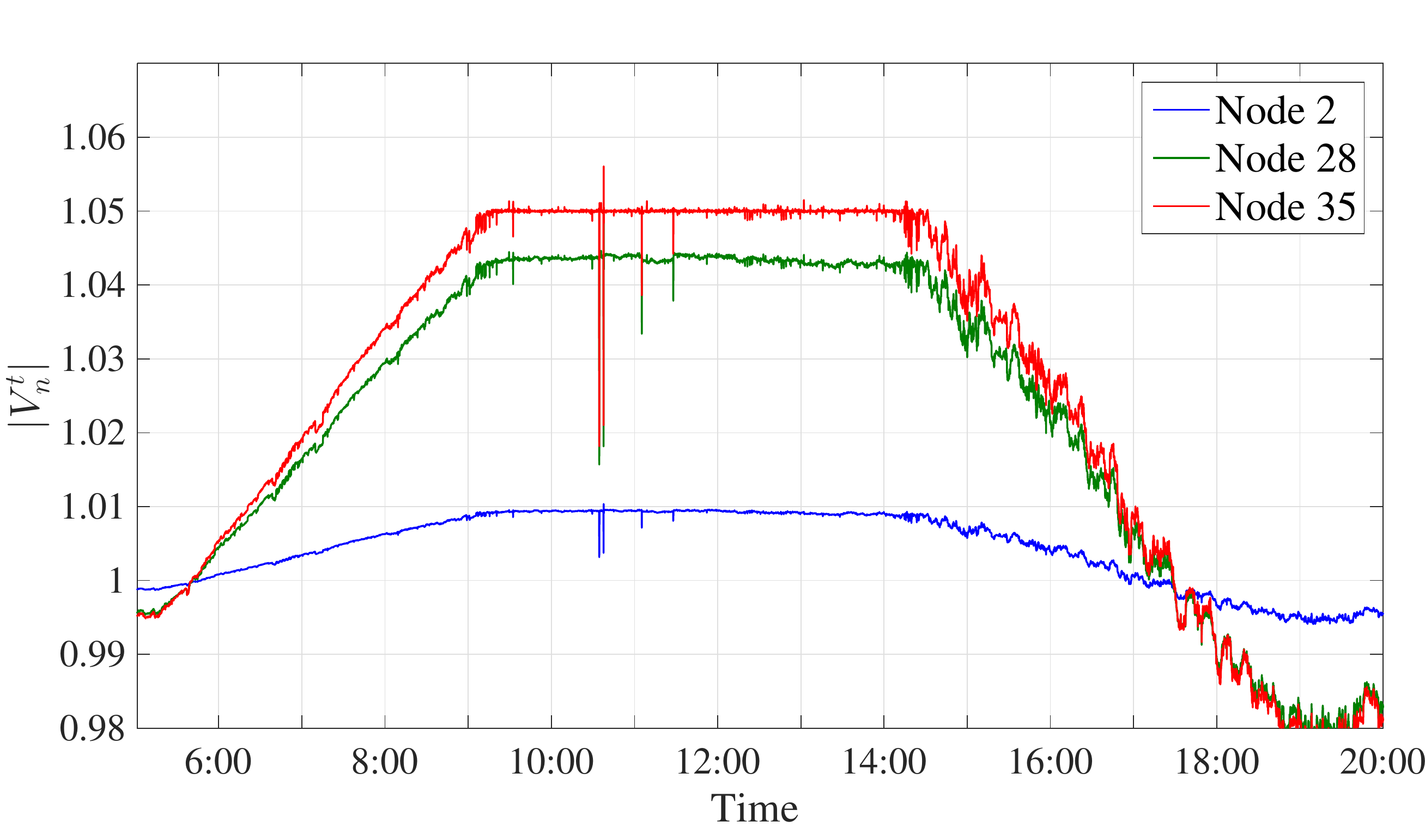} }
\caption{Achieved voltage profile: (a) without control; (b) implementing Volt/Var local control without dead band; and, (c) implementing the proposed controllers.}
\label{F_voltage}
\end{figure}

To this end, consider a modified version of the IEEE 37-node test feeder shown in Figure~\ref{F_feeder}. The modified network is obtained by considering a single-phase equivalent, and by replacing the loads specified in the original dataset with real load data measured from feeders in Anatolia, CA during the week of August 2012~\cite{Bank13}. Particularly, the data have a granularity of $1$ second, and represent the loading of secondary transformers. Line impedances, shunt admittances, as well as active and reactive loads are adopted from the respective dataset. With reference to Figure~\ref{F_feeder}, it is assumed that eighteen  PV systems are located at nodes $4$, $7$, $10$, $13$, $17$, $20$, $22$, $23$, $26$, $28$, $29$, $30$, $31$, $32$, $33$, $34$, $35$, and $36$, and their generation profile is simulated based on the real solar irradiance data available in~\cite{Bank13}. Solar irradiance data have a granularity of $1$ second. The rating of these inverters are $300$ kVA for $i = 3$, $350$ kVA for $i = 15, 16$, and $200$ kVA for the remaining inverters. With this setup, when no actions are taken to prevent overvoltages, one would obtain the voltage profile illustrated in Figure~\ref{F_voltage}(a). To facilitate readability, only three voltage profiles are provided.

\begin{figure}[t] 
\hspace{-.2cm}\includegraphics[width=9.0cm]{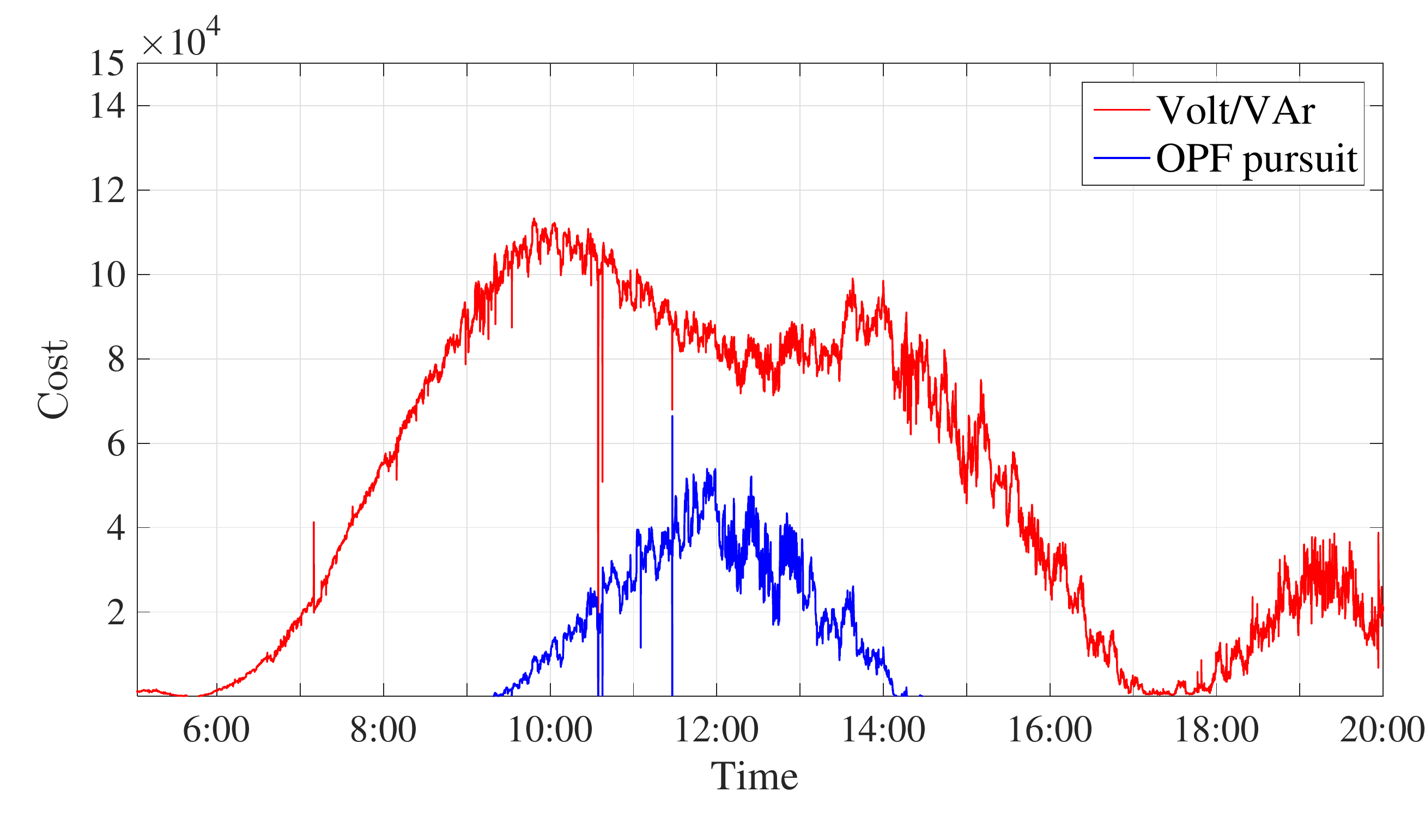}
\caption{Achieved cost $\bar{f}^k(\bu^k) = \sum_{i \in \cG} c_q (Q_i^k)^2 + c_p (P_{\textrm{av},i}^k - P_i^k)^2$.}
\label{F_cost}
\end{figure}

The voltage limits $V_{\mathrm{max}}$ and $V_{\mathrm{min}}$ are set to $1.05$ pu and $0.95$ pu, respectively.  The performance of the proposed scheme is compared against the one of local Volt/VAr control~\cite{Aliprantis13,Zhang13}, one of the control strategies currently tested on the field by a number of DMS vendors and utility companies; particularly, a droop control without deadband~\cite{Aliprantis13,Zhang13} is tested, where inverters set $Q_n^k = 0$ when $|V_n^k| = 1$ pu and linearly increase the reactive power to $Q_n^k = -\sqrt{S_n^2 - (P_{\textrm{av},n}^k)^2}$ when $|V_n^k| \geq 1.05$ pu. The PV-inverters measure the voltage magnitude and update the reactive setpoint every 0.33 seconds. 

For the proposed controllers, the parameters are set as $\nu = 10^{-3}$, $\epsilon = 10^{-4}$, and $\alpha = 0.2$. The stepsize $\alpha$ was selected empirically. The target optimization objective~\eqref{mg-cost2} is set to $\bar{f}^k(\bu^k) = \sum_{i \in \cG} c_q (Q_i^k)^2 + c_p (P_{\textrm{av},i}^k - P_i^k)^2$ in an effort to minimize the amount of real power curtailed and the amount of reactive power injected or absorbed. The coefficients are set to $c_p = 3$ and $c_q = 1$. Iteration of the controllers~\eqref{eq:controllers} is performed every 0.33 seconds. Before describing the obtained voltage profiles, it is prudent to stress that from \emph{Theorem~1} it is evident that the convergence of the controllers \emph{is not} affected by the network size.

Figure~\ref{F_voltage}(b) illustrates the voltage profiles for nodes $2, 28$, and $35$ when Volt/VAr control is implemented. The maximum values of the voltage magnitude are obtained at node $35$. It can be seen that Volt/VAr control enforces voltage regulation, except for the interval between 11:30 and 13:00. In fact, the available reactive power is upper bounded by $(S_n^2 - (P_{\textrm{av},n}^k)^2)^{\frac{1}{2}}$, and this bound decreases with the increasing of $P_{\textrm{av},n}^k$; it follows that in the present test case the inverters do not have sufficient reactive power between 11:30 and 13:00 to enforce voltage regulation. Figure~\ref{F_voltage}(c) shows the voltage profile obtained with the proposed controllers~\eqref{eq:controllers}. It can be seen that the proposed controllers enforce voltage regulation, and a flat voltage profile is obtained from 9:30 to 14:00 [cf. Figure~\ref{F_voltage}(a)]. A flat voltage profile is obtained because in the present test case the controllers minimize the amount of real power curtailed and the amount of reactive power provided; thus, the objective of the controllers is to ensure voltage regulation while minimizing the deviation from the point $[P_{\textrm{av},n}^k, 0]^\sfT$. A few flickers are experienced due to rapid variations of the solar irradiance, but the voltage magnitudes are enforced below the limit within 1-2 seconds.  

Figure~\ref{F_cost} reports the cost achieved by the proposed controllers; that is $\sum_{i \in \cG} c_q (Q_i^k)^2 + c_p (P_{\textrm{av},i}^k - P_i^k)^2$. This is compared against the cost of reactive power provisioning entailed by Volt/VAr control, which is computed as $\sum_{i \in \cG} c_q (Q_i^k)^2$. The advantages of the proposed controllers are evident, as they enable voltage regulation with minimal curtailment of real power as well as reactive power support. Notice that the lower is the amount of reactive power absorbed by the inverters, the lower are the currents on the distribution lines, with the due benefits for the distribution system operators~\cite{Tonkoski11,OID}. It is also worth emphasizing that the cost entailed by Volt/VAr is decreasing during solar-peak hours; as mentioned above, this is because the available reactive power is upper bounded by $(S_n^2 - (P_{\textrm{av},n}^k)^2)^{\frac{1}{2}}$, and this bound decreases with the increasing of $P_{\textrm{av},n}^k$. However, while the cost decreases around 10:00 -- 12:00, the Volt/VAr controllers are not able to ensure voltage regulation.

\begin{figure}[t] 
\hspace{-.2cm}\includegraphics[width=9.0cm]{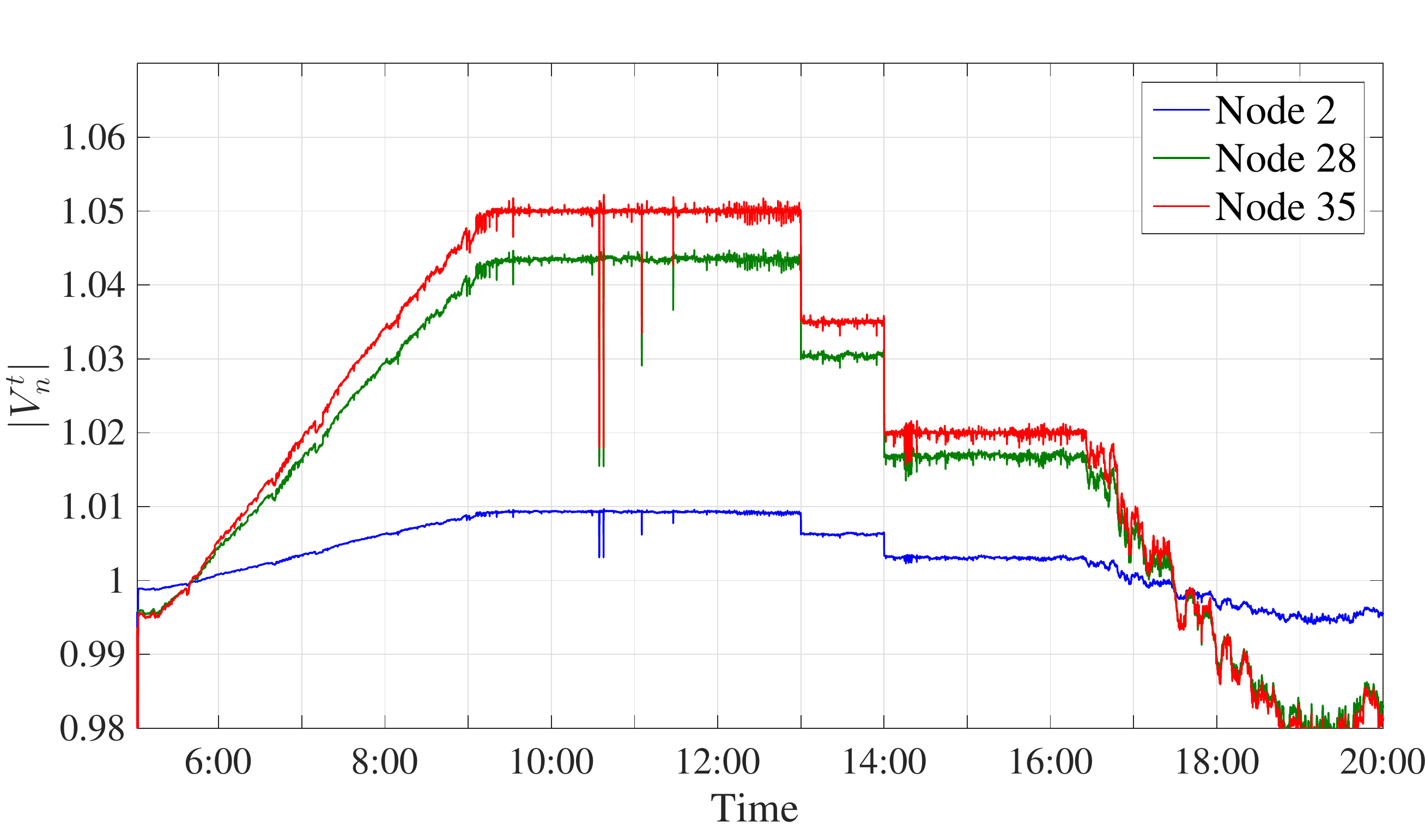}
\caption{Voltage profile achieved by the proposed controllers when  $V^{\mathrm{max}}$ is time-varying. Particularly, $V^{\mathrm{max}}$ is set to: 1.05 pu from 6:00 to 13:00; 1.035 from 13:00 to 14:00; and, 1.02 after 14:00.}
\label{F_V_control_steps}
\vspace{-.5cm}
\end{figure}

Notice that the voltage magnitudes can be forced to flatten on a different value (e.g., $1.045$ pu) by simply adjusting $V^{\mathrm{max}}$. Given the obtained trajectories, it is evident that the proposed  controllers can be utilized to effect also Conservation Voltage Reduction by appropriately changing the values of $V^{\mathrm{min}}$ and $V^{\mathrm{max}}$ in the control loop [cf.~\eqref{eq:updategammaC}--\eqref{eq:updatemuC}]. To test the ability of the proposed controllers to modify the voltage profile in real time in response to changes in $V^{\mathrm{min}}$ and $V^{\mathrm{max}}$, consider the case where the distribution system operator sets the bound $V^{\mathrm{max}}$ to: i) 1.05 pu from 6:00 to 13:00; ii) 1.035 from 13:00 to 14:00; and, iii) 1.02 after 14:00. Figure~\ref{F_V_control_steps} illustrates the voltage profile obtained by the proposed controllers in the present setup. It can be clearly seen that the voltages are quickly regulated within the desired bounds.

\section{Concluding remarks}
\label{sec:conclusions}

This paper addressed the synthesis of feedback controllers that seek RES setpoints corresponding to AC OPF solutions. Appropriate linear approximations of the AC power flow equations were utilized along with primal-dual methods to develop fast-acting low-complexity controllers that can be implemented onto microcontrollers that accompany interfaces of gateways and inverters.  The tracking capabilities of the proposed controllers were analytically established and numerically corroborated. 

\appendix

\subsection{Proof of Theorem~\ref{theorem.inexact}}

For notational simplicity, recall that $\bz^k:= [(\bu^k)^\sfT, (\bgamma^k)^\sfT, (\bmu^k)^\sfT]^\sfT$ collects the primal and dual variables produced by the proposed controllers, and define the following time-varying mapping
\begin{equation*}
\bPhi^k_e: \{\bu^{k}, \bgamma^{k}, \bmu^{k}\} \mapsto \
 \left[\begin{array}{c}
 \nabla_{\bu_1}  \cL_{\nu, \epsilon}^k(\bu, \bgamma, \bmu) |_{\bu_1^{k}, \bgamma^k, \bmu^k} \\
 \vdots \\
 \nabla_{\bu_{N_{\cG}}}  \cL_{\nu, \epsilon}^k(\bu, \bgamma, \bmu) |_{\bu_{N_{\cG}}^{k}, \bgamma^k, \bmu^k}\\
- (V^{\mathrm{min}} -  y_1^k - \epsilon \gamma_1^{k})  \\
\vdots \\
- (V^{\mathrm{min}} -  y_M^k - \epsilon \gamma_M^{k})  \\
 - (y_1^k - V^{\mathrm{max}}  - \epsilon \mu_1^{k}) \\
 \vdots\\
  - (y_M^k - V^{\mathrm{max}}  - \epsilon \mu_M^{k})  
 \end{array}\right] , \hspace{-.1cm}
\end{equation*}
which allows us to rewrite~\eqref{eq:controllers} in the following compact form
\begin{equation}\label{eq:controllers_compact}
\bz^{k+1} = \mathrm{proj}_{\cY^{k} \times \mathbb{R}_{+}^M \times \mathbb{R}_{+}^M}\left\{{\bz}^k - \alpha \bPhi^k_e({\bz}^k)\right\}.
\end{equation}
Consider the norm $\|\bz^k - \bz^{*,k-1}\|_2$, which captures the distance between $\bz^k$ and  the optimal triplet $(\bu^{*,k-1},\bgamma^{*,k-1},\bmu^{*,k-1})$ at time $k-1$ for the min-max problem~\eqref{eq:saddlepoint}. Using~\eqref{eq:controllers_compact}, we can write
\begin{align}
& \|\bz^k - \bz^{*,k-1}\|_2 = \nonumber \\ 
& \left\|\mathrm{proj}_{\cY^{k-1}\times \mathbb{R}_{+}^M \times \mathbb{R}_{+}^M}\left\{{\bz}^{k-1} - \alpha \bPhi^{k-1}_e({\bz}^{k-1})\right\} - \bz^{*,k-1}\right\|_2. \label{dummy_rv1}
\end{align}
By standard optimality conditions, the optimizer is a fixed point of the iterations~\eqref{eq:updateopt_compact}, i.e., $\bz^{*,k-1} = \mathrm{proj}_{{\cY^{k-1}} \times \mathbb{R}_{+}^M \times \mathbb{R}_{+}^M}\left\{\bz^{*,k-1}  - \alpha \bPhi^k(\bz^{*,k-1} )\right\}$. By virtue of this fact,~\eqref{dummy_rv1} can be rewritten as
\begin{multline}\label{dummy_rv2}
\!\!\|\bz^k - \bz^{*,k-1}\|_2 = \left\|\mathrm{proj}_{\cY^{k-1} \times \mathbb{R}_{+}^M \times \mathbb{R}_{+}^M}\left\{{\bz}^{k-1}\!\!\! - \alpha \bPhi^{k-1}_e({\bz}^{k-1})\right\} \right. \\ - \left. \mathrm{proj}_{\cY^{k-1}\times \mathbb{R}_{+}^M \times \mathbb{R}_{+}^M}\left\{\bz^{*,k-1} - \alpha \bPhi^{k-1}(\bz^{*,k-1}) \right\}\right\|_2.
\end{multline} 
We now utilize the non-expansivity property of the projection operator, which yields
\begin{multline}\label{eq.slack1}
\|\bz^{k} - \bz^{*,k-1}\|_2 \leq \|\bz^{k-1} - \alpha \bPhi^{k-1}_e(\bz^{k-1}) \\ - \bz^{*,k-1} + \alpha\bPhi^{k-1}(\bz^{*,k-1})\|_2\ . 
\end{multline}
By construction, observe that
\begin{equation}
\bPhi^{k-1}_e(\bz^{k-1}) -  \bPhi^{k-1}(\bz^{k-1}) = \be^{k-1} 
\end{equation}
where $\be^k := [\mathbf{0}_{2 N_{\cG}}^\sfT, (\be^k_\gamma)^\sfT, (\be^k_\bmu)^\sfT]^\sfT$ is the gradient error. By this definition, we can now expand and bound the right-hand side of~\eqref{eq.slack1} as
\begin{multline}\label{eq.slack2}
\hskip-0.4cm\|\bz^{k-1} \!- \alpha \bPhi^{k-1}(\bz^{k-1}) - \bz^{*,k-1} + \alpha \bPhi^{k-1}(\bz^{*,k-1}) - \alpha\be^{k-1}\|_2 \leq \\ \|\bz^{k-1} - \alpha \bPhi^{k-1}(\bz^{k-1}) - \bz^{*,k-1} + \alpha\bPhi^{k-1}(\bz^{*,k-1})\|_2 + \\ \|\alpha \be^{k-1}\|_2
\end{multline}
where we have used the Triangle inequality.  

We use now Lemma~\ref{lemma-Phi}: first the mapping $\bPhi^k$ is strongly monotone with constant $\eta$, that is 
\begin{multline}\label{eq_dummy_M}
(\bPhi^{k-1}(\bz^{k-1}) - \bPhi^{k-1}( \bz^{*,k-1}))^\sfT(\bz^{k-1} - \bz^{*,k-1}) \geq \\ \eta \|\bz^{k-1} - \bz^{*,k-1}\|_2^2.
\end{multline}
Second, the mapping $\bPhi^k$ is Lipschitz continuous with constant $L_{\nu,\epsilon}$, which implies
\begin{equation}\label{eq_dummy_L}
\|\bPhi^{k-1}(\bz^{k-1}) - \bPhi^{k-1}( \bz^{*,k-1})\|_2^2 \leq L_{\nu,\epsilon}^2 \|\bz^{k-1} - \bz^{*,k-1}\|_2^2.
\end{equation}
By expanding the squared first term in the right-hand side of~\eqref{eq.slack2} and by using the properties~\eqref{eq_dummy_M}-\eqref{eq_dummy_L}, we can write
\begin{multline}\label{eq.slack3}
\|\bz^{k-1} - \alpha \Phi^{k-1}(\bz^{k-1}) - \bz^{*,k-1} + \alpha\Phi^{k-1}(\bz^{*,k-1})\|_2^2  \leq \\
(1 - 2\alpha \eta + \alpha^2 L_{\nu,\epsilon}^2)\|\bz^{k-1} - \bz^{*,k-1}\|_2^2.
\end{multline}
By putting together the results in \eqref{eq.slack1}, \eqref{eq.slack2}, and \eqref{eq.slack3} as well as the bound on the gradient error $\|\be^{k-1}\|_2\leq \sqrt{2} e$ [cf. \emph{Assumption~5}], we have that 
\begin{align}
\|\bz^k - \bz^{*,k-1}\|_2 & \leq \sqrt{{2}}\alpha e + \nonumber \\ &  \sqrt{1 - 2\alpha \eta + \alpha^2 L_{\nu,\epsilon}^2} \|\bz^{k-1} - \bz^{*,k-1}\|_2.
\end{align}
For simplicity, let $\rho(\alpha) : = \sqrt{1 - 2\alpha \eta + \alpha^2 L_{\nu,\epsilon}^2}$. Thus, it follows that
\begin{align}\label{dummy1}
\|\bz^k - \bz^{*,k-1}\|_2 \leq \rho(\alpha) \|\bz^{k-1} - \bz^{*,k-1}\|_2 + \sqrt{2}\alpha e.
\end{align}
We now consider the distance between the controller output $\bz^k$ with the current optimizer of the min-max problem~\eqref{eq:saddlepoint}, i.e., $\|\bz^k - \bz^{*,k}\|_2$. This quantity can be bounded by using \emph{Assumptions~3--4} on the variability of primal and dual optimizers. Particularly, by using the Triangle inequality and the relation~\eqref{dummy1}, it follows that
\begin{align}\label{dummy2}
\|\bz^k - \bz^{*,k}\|_2 & = \|\bz^k - \bz^{*,k} - \bz^{*,k-1} + \bz^{*,k-1}\|_2 \nonumber \\ & \leq \|\bz^k - \bz^{*,k-1}\|_2 + \sigma_{\bz} \nonumber \\
& \leq \rho(\alpha) \|\bz^{k-1} - \bz^{*,k-1}\|_2 + \sqrt{2}\alpha e + \sigma_{\bz}.
\end{align}
If $\rho(\alpha)< 1$, then~\eqref{dummy2} represents a contraction, and via the geometric series sum formula we can write 
\begin{align}\label{dummy3}
\|\bz^k - \bz^{*,k}\|_2 & \leq [\rho(\alpha)]^k \|\bz^0 - \bz^{*,0}\|_2 + \nonumber \\
& \quad \frac{1 - [\rho(\alpha)]^k}{1 - \rho(\alpha)} \Big[\sqrt{2}\alpha e + \sigma_{\bz}\Big]. \nonumber 
\end{align}
The relation above describes a Q-linear convergence of $\|\bz^k - \bz^{*,k}\|_2$ to a neighborhood of $0$, with asymptotic error bound given by
\begin{equation}
\limsup_{k\to \infty} \|\bz^k - \bz^{*,k}\|_2 = \frac{1}{1 - \rho(\alpha)} \Big[\sqrt{2}\alpha e + \sigma_{\bz}\Big],
\end{equation}
which completes the proof.

\subsection{Setpoint update}
\label{sec:setpointclosedform}

The setpoint update~\eqref{eq:updatevC} affords a closed-form solution for a variety of RESs and other controllable devices. For notational simplicity, let $\hat{\bu}_n^{k} = [\hat{P}_n^{k}, \hat{Q}_n^{k}]^\sfT$ be the unprojected point, where $\hat{P}_n^{k}$ and $\hat{Q}_n^{k}$ are the unprojected values for the real and reactive powers, respectively; that is, 
\begin{align} 
\hat{\bu}_n^{k} :=  \bu_n^{k-1}  - \alpha \nabla_{\bu_n}  \cL_{\nu, \epsilon}^{k-1}(\bu, \bgamma, \bmu) |_{\bu_n^{k-1}, \bgamma^{k-1},\bmu^{k-1}} \, .
\end{align}
Clearly, one has that $\bu_n^{k}  =  \mathrm{proj}_{\cY_n^{k-1}}\{\hat{\bu}_n^{k} \}$. In the following, expressions for $\bu_n^{k}$ are reported for different choices of the set $\cY_n^{k-1}$. 

\noindent \emph{Real power-only control}: in this case, the set $\cY_n^{k-1}$ boils down to $\cY_n^{k-1}  =  \left\{(P_n, Q_n) \hspace{-.1cm} :  0  \leq {P}_{n}  \leq  P_{\textrm{av},n}^{k-1}, Q_n = 0 \right\}$. This set is typical in inverter-interfaced RESs adopting real power curtailment-only strategies~\cite{Tonkoski11,vonAppen14}, where $P_{\textrm{av},n}^{k-1}$ represents the maximum power point for a RES; it also represents conventional generation unit operating at unity power factor, where $P_{\textrm{av},n}^{k-1}$ is the maximum generation. In this case,~\eqref{eq:updatevC} can be simplified as follows:
\begin{subequations}
\begin{align}
P_n^{k} & =  \max\{0, \min\{\hat{P}_n^{k}, P_{\textrm{av},n}^{k-1} \}   \}\\
Q_n^{k} & = 0 \, .
\end{align}
\end{subequations}

\noindent \emph{Reactive power-only control}: For RES with reactive power-only control capability, the set of possible operating points is given by $\cY_n^{k-1} = \{(P_n, Q_n)  \hspace{-.1cm} : P_n = P_{\textrm{av},n}^{k-1},  |Q_{n}|  \leq  (S_{n}^2 - (P_{\textrm{av},n}^{k-1})^2)^{\frac{1}{2}} \}$~\cite{Aliprantis13,Farivar12}. In this case,~\eqref{eq:updatevC} boils down to:
\begin{subequations}
\begin{align}
P_n^{k} & = P_{\textrm{av},n}^{k-1} \\
Q_n^{k} & = \mathrm{sign}(\hat{Q}_n^{k}) \min\{|\hat{Q}_n^{k}|,  (S_{n}^2 - (P_{\textrm{av},n}^{k-1})^2)^{\frac{1}{2}}\}
\end{align}
\end{subequations}
where $\mathrm{sign}(x) = -1$ when $x < 0$ and $\mathrm{sign}(x) = 1$ when $x > 0$.

\begin{figure}[t] 
\begin{center}
\includegraphics[width=5.5cm]{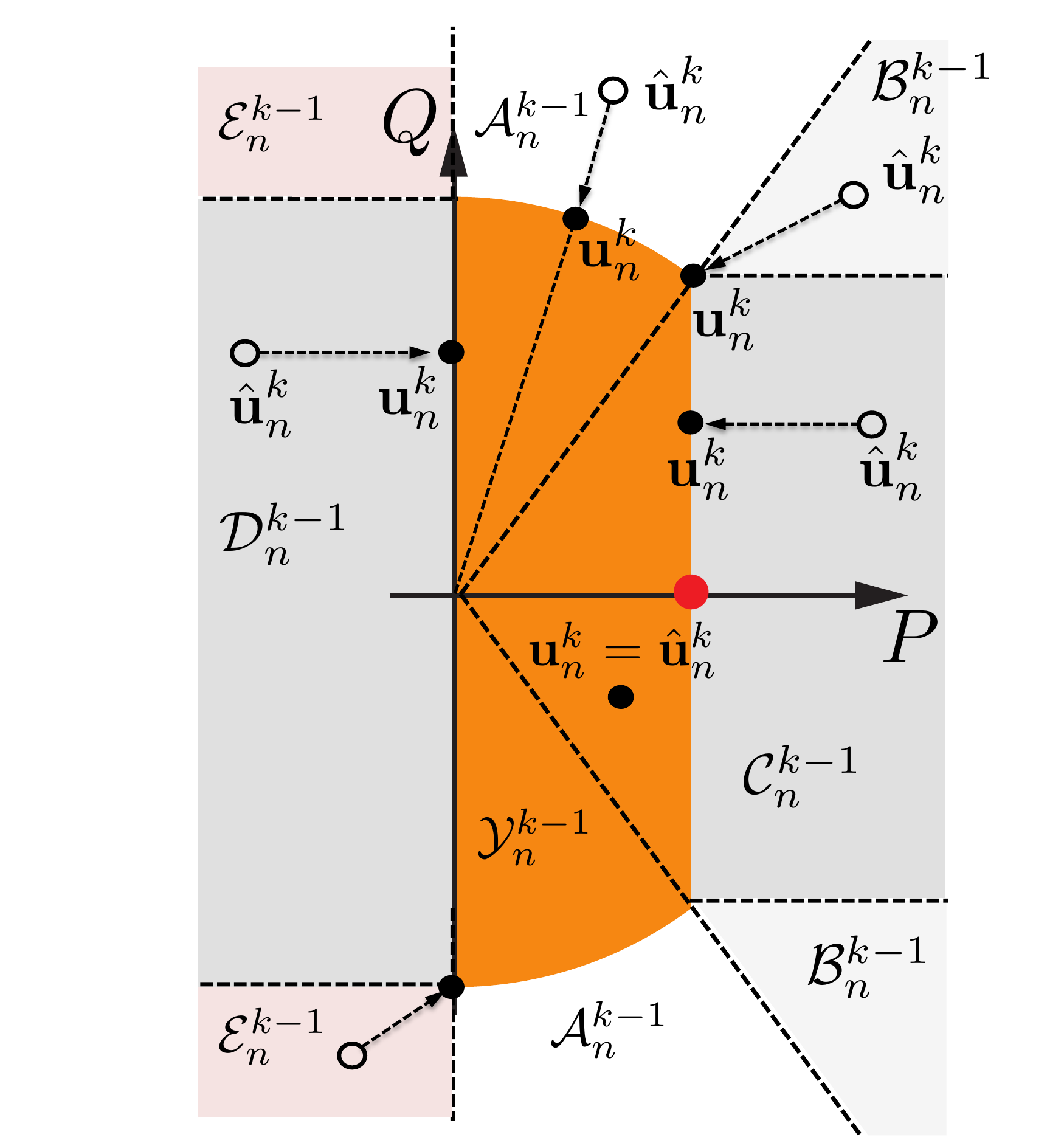}
\end{center}
\vspace{-.5cm}
\caption{Projection onto set~\eqref{mg-PV}. The red dot corresponds to the point $[P_{\textrm{av},i}^{k-1}, 0]^\sfT$.}
\label{F_proj}
\vspace{-.5cm}
\end{figure}

\noindent \emph{Joint real and reactive power control}: Consider now the more general setting where an inverter can control both real and reactive output powers; particularly, given the inverter rating $S_n$ and the current available real power $P_{\textrm{av},n}^{k-1}$, consider the set $ \cY_n^{k-1}  =  \left\{({P}_{n}, {Q}_{n} ) \hspace{-.1cm} :  0 \leq {P}_{n}  \leq  P_{\textrm{av},n}^{k-1}, ({Q}_{n})^2  \leq  S_{n}^2 - ({P}_{n})^2 \right\} $ in~\eqref{mg-PV}. With reference to Figure~\ref{F_proj}, the setpoints $\bu_n^{k}$ can be obtained from the unprojected point $\hat{\bu}_n^{k}$ as summarized next:
\begin{align}
\bu_n^{k} = 
\left\{
\begin{array}{ll}
\hat{\bu}_n^{k} \hspace{2.25cm} \textrm{, ~if~} \hat{\bu}_n^{k} \in \cY_n^{k-1} \\
\hat{\bu}_n^{k} \frac{S_n}{\|\hat{\bu}_n^{k} \|} \hspace{1.55cm} \textrm{, ~if~} \hat{\bu}_n^{k} \in \cA_n^{k-1} \\
\left[ P_{\textrm{av},n}^{k-1}, \mathrm{sign}(\hat{Q}_n^{k}) (S_{n}^2 - (P_{\textrm{av},n}^{k-1})^2)^{\frac{1}{2}} \right]^\sfT ,  \\
\hspace{2.85cm} \textrm{~if~} \hat{\bu}_n^{k} \in \cB_n^{k-1} \\
\left[ P_{\textrm{av},n}^{k-1}, \hat{Q}_n^{k} \right]^\sfT \hspace{.75cm} \textrm{, ~if~} \hat{\bu}_n^{k} \in \cC_n^{k-1}  \\
\left[ 0, \hat{Q}_n^{k} \right]^\sfT \hspace{1.3cm} \textrm{, ~if~} \hat{\bu}_n^{k} \in \cD_n^{k-1}  \\
\left[ 0, \mathrm{sign}(\hat{Q}_n^{k}) S_{n} \right]^\sfT  \textrm{, ~if~} \hat{\bu}_n^{k} \in \cE_n^{k-1}  \\
\end{array}
\right.
\end{align}
where the regions $\cA_n^{k-1}, \cB_n^{k-1}, \cC_n^{k-1}$, $\cD_n^{k-1}$, and $\cE_n^{k-1}$ can be readily obtained from  $S_n$ and $P_{\textrm{av},n}^{k-1}$.

It is also worth pointing out that closed-form expressions can be found when $\cY_n^{k-1}$ models the operating regions of e.g., diesel generators and controllable loads with variable speed drives.

\bibliographystyle{IEEEtran}
\bibliography{biblio.bib}

\end{document}